\RequirePackage{etoolbox}
\csdef{input@path}{%
 {sty/}
 {img/}
}%

\documentclass[ba]{imsart}
\pubyear{0000}
\volume{00}
\issue{0}
\doi{0000}
\firstpage{1}
\lastpage{1}
\usepackage[T1]{fontenc}
\usepackage[utf8x]{inputenc}
\usepackage{amsmath}
\usepackage{amssymb}
\usepackage{graphicx} 

\usepackage{amsthm}
\usepackage{color}
\usepackage{float}
\usepackage{multirow}
\graphicspath{{./figures/}}

\newtheorem{theorem}{Theorem}
\newtheorem{lemma}{Lemma}

\newcommand{\R}{\mathbb{R}}

\newcommand{\F}{\ensuremath{\mathcal{F}}}
\newcommand{\eps}{\varepsilon}

\renewcommand{\P}{\ensuremath{\mathcal{P}}}

\usepackage{lmodern}

\DeclareMathAlphabet{\mathpzc}{OT1}{pzc}{m}{it}

\newcommand{\I}{\mathbb{I}}

\newcommand{\E}{\mathrm{E}}
\newcommand{\N}{\mathbb{N}}

\usepackage{float}

\usepackage{nameref}
\usepackage{natbib}

\usepackage{stmaryrd}
\usepackage{graphicx}

\makeatletter
\newcommand{\fixed@sra}{$\vrule height 2\fontdimen22\textfont2 width 0pt\shortrightarrow$}
\newcommand{\shortarrow}[1]{%
  \mathrel{\text{\rotatebox[origin=c]{\numexpr#1*45}{\fixed@sra}}}
}
\makeatother

\newcommand{\M}{\mathcal{M}}
\newcommand{\FO}{\mathcal{F}_{+}}
\newcommand{\Fm}{\mathcal{F}_{\shortarrow{7}}}

\begin{document}

\begin{frontmatter}
\title{Testing un-separated hypotheses by estimating a distance}
\runtitle{Bayes test}

\begin{aug}
\author{\fnms{Jean-Bernard} \snm{Salomond}\thanksref{addr1,t1}\ead[label=e1]{jean-bernard.salomond@u-pec.fr}}

\runauthor{J-B Salomond}

\address[addr1]{Universit\'e Paris-Est, Laboratoire d'Analyse et de Math\'ematiques Appliqu\'ees (UMR 8050), UPEM, UPEC, CNRS, F-94010, Cr\'eteil, France}
    \printead{e1}

\end{aug}

\begin{abstract}
In this paper we propose a Bayesian answer to testing problems when the hypotheses are not well separated. The idea of the method is to study the posterior distribution of a discrepancy measure between the parameter and the model we want to test for. This is shown to be equivalent to a modification of the testing loss. An advantage of this approach is that it can easily be adapted to complex hypotheses testing which are in general difficult to test for. Asymptotic properties of the test can be derived from the asymptotic behaviour of the posterior distribution of the discrepancy measure, and gives insight on possible calibrations.
In addition one can derive separation rates for testing, which ensure the asymptotic frequentist optimality of our procedures. 
\end{abstract}


\end{frontmatter}

\section{Introduction}
\label{sec:intro}


{
Bayesian hypothesis testing, although widely studied in the literature, is still subject to controversy \citep[see][to name a few]{Jeffreys1939, Bernardo1980, BergerSellke,gelman2008objections}. 
In particular, a lot of efforts have been puted on reconciling Bayesian and frequentist testing procedures as in \cite{BergerSellke}, \cite{BergerUnified} or \cite{BergerPrecise87}. In this paper, we focus on the specific case of two-hypotheses testing although we believe that the ideas developed here are more general; more precisely, we consider testing problems of the form: 
\begin{equation}
H_0 : \theta \in \M_0, ~\text{versus} ~ H_1: \theta \in \M_1,
\label{eq:test:generl}
\end{equation}
where $\M_0$ and $\M_1$ are not well separated, i.e. $\bar{\M}_0 \cap \bar{\M}_1 \neq \emptyset$ where $\bar{\mathcal{F}}$ stands for the closure of $\mathcal{F}$. 
When considering prediction, it is now well known that standard Bayesian methods such as BIC have a tendency to favour the simpler model, even when the more complex one gives better predictions, as shown in \cite{RSSB:RSSB1025}. This phenomenon also occurs in a testing or model selection setting when hypotheses are nested, and induces a loss of power for the Bayesian test near the null. 
In our view, one reason for this lack of efficiency of standard Bayesian testing approaches, such as the Bayes Factor or the comparison of posterior probabilities, comes from the fact that parameters that are close to the boundary between both hypotheses can be approximated from both sides. Thus, depending on the prior distribution on both the null and the alternative, some inconsistency may occur. This phenomenon is shown on some examples in section \ref{sec:boundary}. This loss of power of Bayesian testing procedures, induced by the prior is troublesome as it is difficult to control for and strongly depends on the prior distribution. 
Finding good prior distributions for testing has been a subject of high interest in the recent years. In particular, \cite{JohnsonRossell2010} \citep[or the actualised version of their ideas developed in][]{doi:10.1080/01621459.2015.1130634} consider a similar case of un-separated hypotheses. Their idea is to enforce separation through the prior distribution using \emph{non local} priors. As exposed in \cite{rousseau2010moment}, this approach can be viewed as a modification of the loss used for testing. It appears on simple examples studied in section \ref{sec:boundary} that imposing such a penalty can make it more difficult to detect parameters near the boundary between hypotheses (see section \ref{sec:pointnull} for instance).
The problem of finding a good prior distribution for testing has also been tackled by \cite{johnson2013}, where the author introduced uniformly most powerful Bayesian test. The author proposes to calibrate the method by maximizing the probability that the Bayes Factor exceeds a certain threshold under the alternative. However, the proposed method seems difficult to extend outside exponential models. 
In this paper we propose a novel approach to the problem of testing un-separated hypotheses, based on the evaluation of a \emph{discrepancy} between the parameter $\theta$ and the hypothesis at hand. 
A great advantage of the approach is that it is in general easy to use in practice and it generalizes directly to nonparametric hypotheses testing. 
Let $D(\theta,\M_0)$ be a \emph{discrepancy measure} between $\theta$ and $\M_0$. Following the frequentist approach to testing, our idea is to associate $\theta$ to $\M_0$ if $D(\theta,\M_0)$ is bellow a certain threshold $\tau$. This idea of choosing the model closer to the parameter for a certain metric is quite general and we believe that it could be applied in a wide variety of settings. In this paper, we might only focus on the simpler problem of two hypotheses testing.  

Although not aiming at the same problem, this approach is similar to the idea of approximating precise hypotheses by point null hypotheses as studied in \cite{BergerPrecise87}, which can be re-interpreted as a use of non-local prior as argued in \citet{JohnsonRossell2010}. This approximation of hypotheses where latter studied in \cite{verdinelli1998bayesian} and \cite{Rousseau2007}. More specifically, in the latter the author proposes a generalization of the $0-1$ loss function from which a Bayesian test is derived, and which induces a separation of the hypotheses. Following \cite{Rousseau2007}, we consider the following loss function 
\begin{equation}
L(\theta,\delta) = 
\begin{cases}
0& \text{ if } \delta = \I_{D(\theta,\M_0) \leq \tau} \\ 
\gamma_0 & \text{ if } \delta = 1 \text{ and }  D(\theta,\M_0) \leq \tau \\ 
\gamma_1 & \text{ if } \delta = 0 \text{ and }  D(\theta,\M_0) > \tau \\ 
\end{cases},
\label{eq:loss:modified}
\end{equation}   
where the parameters $\gamma_0$ and $\gamma_1$ have the same interpretation as for the standard weighted $0-1$ loss \citep[see][]{Robert2007} in terms of price of misclassification error. A default choice is to take $\gamma_0 = \gamma_1$. 
This modification of the loss function can also be viewed as a relaxation of the hypotheses 
\begin{equation}
H_0^a : D(\theta,\M_0) \leq \tau, ~ H_1^a:  D(\theta,\M_0) > \tau,
\label{eq:test:relaxed}
\end{equation}
For a fixed threshold $\tau$, the same idea was applied in \cite{Dunson01122008} and  \cite{wangDunson11} for testing equality in distribution against stochastic ordering. From a decision theoretic point of view, this loss is relevant since it indicates that we do not {pay} for misclassified parameters that lie in a region in which we cannot differentiate the null and the alternative. 
In addition, as argued in \cite{BergerPrecise87}, one is in general not so interested in knowing if $\theta$ belongs to $\M_0$ but rather if $\theta \in \M_0$ is reasonable approximation. 
From a more practical point of view, this approach gives a method for constructing Bayesian tests that separate well the hypotheses in a wide variety of contexts, including complex alternatives such as nonparametric models for example. Deriving the Bayesian answer to \eqref{eq:test:relaxed} can also lead to simpler procedures. The Bayesian estimate associated with a prior $\Pi$ on the parameter set $\M = \M_0 \cup \M_1$, the loss \eqref{eq:loss:modified} and data $Y^n$, is given by 
\begin{equation}
\delta_n^\pi(\tau) = 
\begin{cases}
0& \text{ if } \Pi\left\{ D(\theta,\M_0)  \leq \tau | Y^n \right\} \geq \frac{\gamma_0}{\gamma_0 + \gamma_1} \\ 
1& \text{ otherwise}
\end{cases}  
\label{eq:intro:decision:rule}
\end{equation} 
where $\Pi(\cdot|Y^n)$ denote the posterior measure of the parameter $\theta$ given the observations $Y^n$.
From this last equation, we see that the behaviour of a test based on our modified $0-1$ loss is driven by the behaviour of $D(\theta,\M_0)$. This will prove particularly useful when testing complex or nonparametric versus complex or nonparametric hypothesis which is known to be a difficult case to handle and has not received much attention in the Bayesian literature. In addition even for simpler model, the behaviour of $D(\theta,\M_0)$ may also be easy to study for a wide variety of priors, as shown in section \ref{sec:parametric} for instance. 
From this formulation, we see that prior distributions that induce a good behaviour for $D(\theta,\M_0)$ in terms of concentration properties will also be good candidates for testing with this approach. Note that such priors may differ from the one that leads to good properties for estimating $\theta$, as shown for example in section \ref{sec:WN}. Note also that to compute the Bayesian test with formulation \eqref{eq:intro:decision:rule}, we only have to sample under the posterior. 
This thus gives leads to tackle two of the main difficulties in Bayesian testing when studying the Bayes Factor: choosing a appropriate prior and computing the marginal distribution.

Once the discrepancy measure is chosen, the remaining problem is calibrating the threshold $\tau$. In an informative context where one has prior knowledge on acceptable discrepancy from $\M_0$, $\tau$ can be calibrated subjectively. However, such a prior knowledge may not be available. We thus propose a calibration of $\tau$ based on asymptotic arguments.
Heuristically, one would like to find a threshold $\tau$ that minimizes the testing error. \cite{johnson2013} proposed a similar idea for constructing uniformly most powerful Bayesian tests where he proposes to chose a prior for testing that maximizes the probability that the Bayes Factor exceed a certain threshold for all $\theta \in \M_1$. 
In our case, in general, minimizing the testing error might not be possible even for some simple models. We thus propose a calibration method based on the asymptotic control of the type I and type II errors. More precisely we chose $\tau = \tau_n$ to be the smallest sequence such that
\begin{equation} 
\sup_{\theta \in \M_0} \E_\theta^n \{\delta^\pi_n(\tau_n)\} = o(1) \label{eq:consistencyH0},
\end{equation}
where $\E_\theta^n$ denote the expectation with respect to $Y^n \sim P_\theta$. Given the formulation of the test \eqref{eq:intro:decision:rule} finding such a calibration will only requires a control of the asymptotic behaviour of $D(\theta,\M_0)$ under the posterior. 
We then study for which sequence of $\rho_n$ we have
\begin{equation}
\sup_{\theta \in \M_1, d(\theta,\M_0)> \rho_n} \E_\theta^n \{1-\delta^\pi_n(\tau_n)\} = o(1) \label{eq:consistencyH1}.
\end{equation}

The sequence $\rho_n$ is thus an upper bound on the separation rate of the test \citep[see][]{lepski2000asymptotically}. Separation rates indicates how close a parameter from $\M_1$ can be to $\M_0$ and sill be detected by the test. Although separation rates have been widely studied in the frequentist literature, to the author's best knowledge, the only related result in the Bayesian literature has been proposed in \cite{doi:10.1080/01621459.2015.1130634}. Note that if the test separates both hypotheses at the best possible rate in the minimax sense, it indicates that the decision rule $\delta_n^\pi(\tau)$ although being a Bayesian answer to the relaxed testing problem \eqref{eq:test:relaxed}, is also an asymptotically optimal frequentist answer for the original testing problem \eqref{eq:test:generl}. This indicates that such a test can catch up with frequentist methods for detecting parameters close to the boundary between the hypotheses. This is to the best of our knowledge a new result for Bayesian test. A counterpart will be of course a loss in parsimony enforcement.  
In the remainder of the paper, we study on two examples the problems that can occurs at or close to the boundary between hypotheses. We then propose a general calibration for $\tau_n$ for some usual testing problems and show that our method achieve the minimax separation rates in these cases. On the last sections, we compare our approach to existing ones for a non-parametric test.

\section{Boundary problems} 
\label{sec:boundary}
In this section we illustrate on simple examples the problems faced by the \emph{non-local prior approach to testing} proposed by \cite{JohnsonRossell2010} and further developped in \cite{doi:10.1080/01621459.2015.1130634} and the \emph{standard priors} when the parameter is at, or near the boundary between the null and the alternative. 

\subsection{Point null hypotheses} 
\label{sec:pointnull}
Consider the following data generating process $X^n \sim \mathcal{N}(\theta,1/\sqrt{ n})$, and the test $H_0: ~\theta = 0$ versus $H_1:~\theta \neq 0$. To compute the standard Bayes Factor for this problem, 
define $\M_0 = \{ 0 \}$ and $\M_1 = \R \backslash \{0\}$, and let the prior distribution $\pi$ on $\theta$ be 
$
\pi : \theta \sim \mathcal{N}(0,\sigma^2) ,
$ 
and chose equal prior weights on both hypotheses. 
We can easilly derive the usual Bayes-Factor for this problem and get 
$$
B_{0,1}(X^n) = {\int_{\R} \pi(\theta) e^{-\frac{n}{2} (X^n - \theta)^2} d\theta \over e^{-\frac{n}{2} (X^n)^2} },
$$
and compare it to $1$. Comparing the Bayes Factor with the fixed threshold $c = 1$ is equivalent to comparing the posterior mass of $\M_0$ with $1/2$. 
For the non local prior we use the method of moment proposed in \cite{doi:10.1080/01621459.2015.1130634} with parameter fixed as proposed in their paper, i.e. 
$$
\pi_1^M (\theta) = \frac{\theta^2}{\tau} \pi(\theta/\tau).
$$ 
The form of the proposed prior is displayed in Figure \ref{fig:pi1m}. We easily derive the Bayes-Factor associated with this prior 
$$
B_{0,1}^M(X^n) = {\int_{\R} \pi_1^M(\theta) e^{-\frac{n}{2} (X^n - \theta)^2} d\theta \over e^{-\frac{n}{2} (X^n)^2} }. 
$$
Here again we shall compare it to 1. 
\begin{figure}[h]
\begin{center}
\includegraphics[scale=0.5]{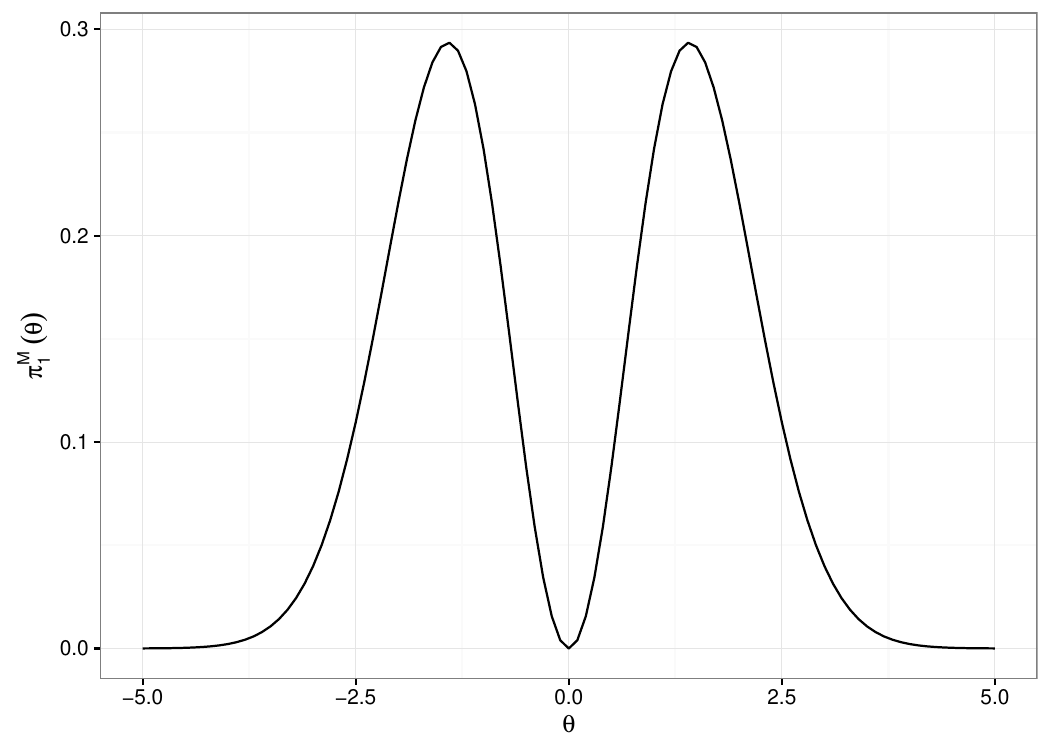}
\caption{Prior on the alternative constructed with the Method of Moment from \cite{doi:10.1080/01621459.2015.1130634}. The hyperparameter $\tau$ is fixed to 0.358.}
\label{fig:pi1m}
\end{center}
\end{figure}
For the method proposed in this paper, we chose as a discrepancy measure $D(\theta,\M_0) = |\theta|$. We now have to calibrate $\tau_n$ such that the test satisfies \eqref{eq:consistencyH0}-\eqref{eq:consistencyH1}. We shall see in the following Theorem \ref{thm:param} that in this case, choosing $\tau_n = u_n n^{-1/2}$ for any $u_n \to + \infty$ will ensure consistency. To calibrate $u_n$, note that we have 
\begin{align}
\pi(D(\theta,\M_0)> \tau_n | X^n) = 1-\Phi\left(\frac{\tau_n - m_x}{\sigma_x}\right) + \Phi\left(\frac{-\tau_n - m_x}{\sigma_x}\right)
\end{align}
where $\sigma_x^2 = (n+\sigma^{-2})^{-1}$ and $m_x = nX^n \sigma^2_x$ are the posterior mean and the posterior variance respectively and $\Phi$ is the cumulative distribution function of a standard Gaussian. To get low type I error while not deteriorating the separation rate we choose $u_n = \max(\Phi^{-1}(0.05),\log(\log(n)))$. We run all three methods on simulated data generated for three different parameters $\theta_0$, namely $\sqrt{2 \log(n)/n}$, $\sqrt{\log(n)/n}$ and $0$. The first two parameters are getting closer and closer to the boundary between hypotheses as the number of observations grows while the third is in $\M_0$. We observe that even when the parameter is at a reasonable distance from $\M_0$, the non-local prior seems to penalize too much, and thus will contract on the simpler model, while the other approaches do detect the parameter as non-zero. When the parameter is at a distance $\sqrt{\log(n)/n}$ then the usual Bayes Factor does not clearly detect the parameter has non null, while the proposed method asymptotically does. The price to pay is a slower decay of the type I error of the order of $\log(n)$ to be compared to a exponential decay for the Bayes Factor.
\begin{figure}[h]
\includegraphics[width = \textwidth]{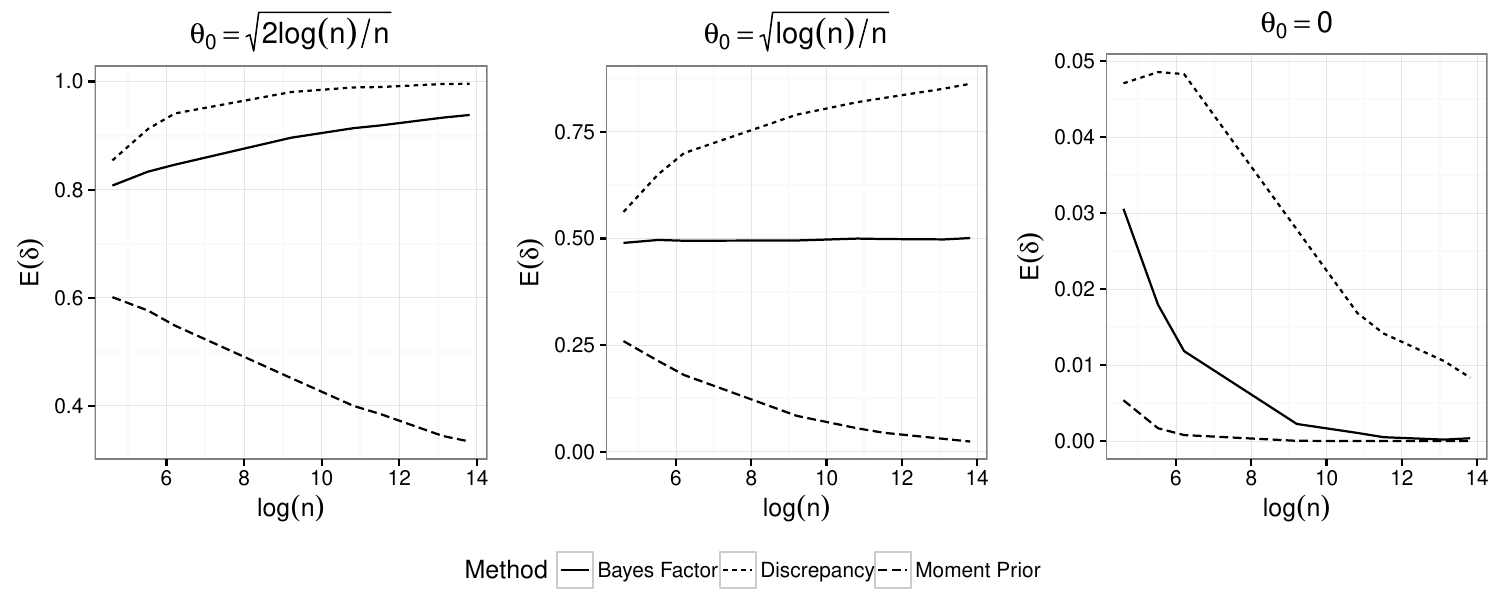}
\begin{center}\caption{Proportion of test that classifies the parameter as non null for $N = 5\times 10^4$ replications of the test. The Bayes Factor obtained with Gaussian and non local priors are compared to 1. For the discrepancy method $\tau_n = \max[1.96,\log(\log(n))]/\sqrt{n}$.}
\label{fig:BFs}
\end{center}
\end{figure}

\subsection{Un-separated hypotheses}

We now consider a case where both the null and the alternative have similar sizes. Using the same setting as before we now test for $H_0: \theta \leq 0$ versus $H_1: \theta > 0$, and thus $\M_0 = (-\infty, 0]$ and $\M_1 = (0,+\infty)$, using the same prior $\pi$ as before. To compute the usual Bayes Factor, we thus have the following on $\M_0$ and $\M_1$ respectively: 
$$
\pi_0(\theta) = 2\pi(\theta) \I_{\theta \in \M_0}, ~ \pi_1(\theta) = 2\pi(\theta) \I_{\theta \in \M_1}.
$$
We can compute Bayes Factor 
$$
B_{0,1}(X^n) = {\int_{\M_0} \pi(\theta | X^n) d\theta  \over \int_{\M_1} \pi(\theta | X^n) d \theta} = \frac{\Phi\left( - m_x /\sigma_x\right)}{\Phi(m_x/\sigma_x)}.
$$
From this formulation, we see that the Bayse Factor $B_{0,1}$ will not detect the parameters $\theta = 0$ (which is at the boundary between $\M_0$ and $\M_1$) as belonging to $\M_0$, leading to poor frequentist performances of such a test in this case. To compare this approach to the non-local prior method
we construct a prior $\pi^M_1$ on the alternative using the method of moment described in \cite{doi:10.1080/01621459.2015.1130634} that will enforce a separation of the hypotheses. We consider the following modification of the prior 
$
\pi^M_1(\theta) = \frac{\theta^2}{\tau} \pi_1(\theta/\tau). 
$
A plot of this prior is given in Figure \ref{fig:pi1mpos}. We can then compute the Bayes Factor 
$$
B_{0,1}^M(X^n) = {\int_{\M_0} \pi_0(\theta) e^{-n(X^n - \theta)^2/2 } d\theta  \over \int_{\M_1} \pi_1^M(\theta) e^{-n(X^n - \theta)^2/2 } d \theta} .
$$
One can easily compute the marginals using simple Monte-Carlo integration. Here again we will compare the Bayes Factor $B_{0,1}^M$ with the fixed threshold $1$. 
\begin{figure}[h]
\begin{center}
\includegraphics[scale=0.5]{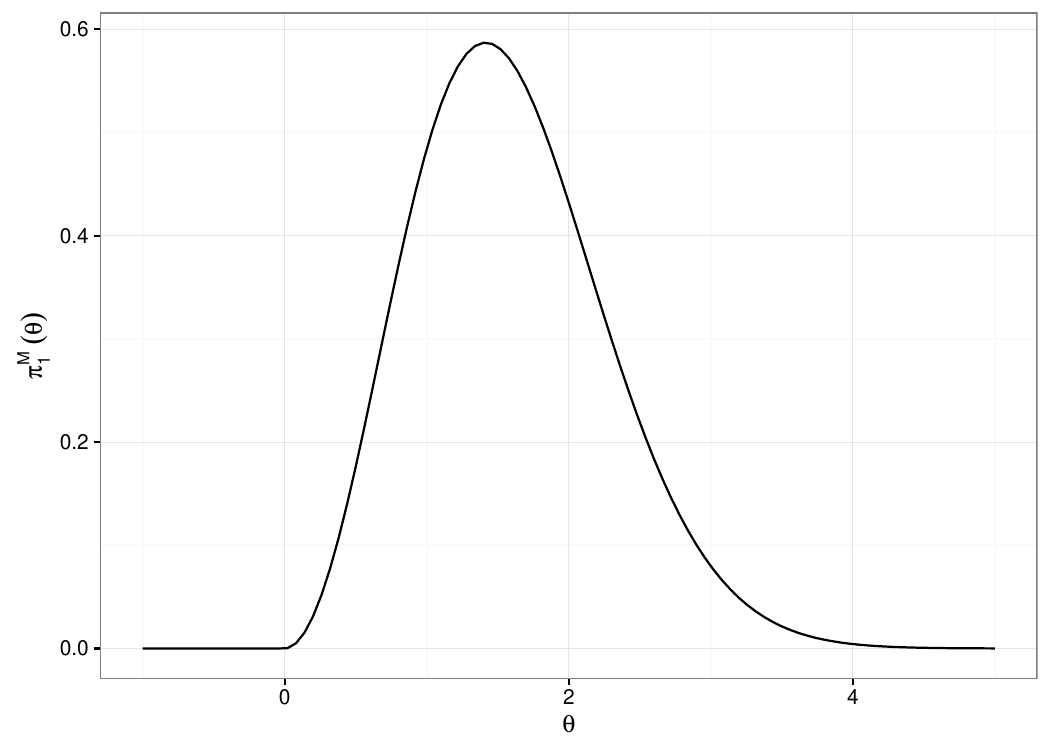}
\caption{Prior on the alternative constructed with the Method of Moment from \cite{doi:10.1080/01621459.2015.1130634}. The hyperparameter $\tau$ is fixed to 0.358.}
\label{fig:pi1mpos}
\end{center}
\end{figure}
In order to compare these approaches with the one proposed in this paper, we first need to find a discrepancy measure $D$ and calibrate the threshold $\tau$. We choose $D(\theta,\M_0) = \min(\theta,0)$. Using a simple standard Gaussian prior, we can easily calibrate the threshold $\tau_n$ using the same approach as before. We have that for all sequence $u_n$ that goes to infinity as slowly as needed, $\tau_n = C u_n n^{-1/2}$ leads to a separation rate $\rho_n \leq 2 \tau_n$. We now calibrate the constant $C$ and the sequence $u_n$ based on heuristics.
Again, denoting $\hat{\sigma}^2_x = (n + 1/\sigma^2)^{-1}$ and $m_x = n X^n \hat{\sigma}^2_x$ the posterior mean and posterior variance respectively, we have %
$$
\Pi(D(\theta,\M_0) \geq \tau_n | X^n) = 1-  \Phi\left( {\tau_n - \hat{m}_x \over \hat{\sigma}_x} \right),
$$
We then choose again $u_n = \max(\log(\log(n)), \Phi^{-1}(0.05))$ which insure consistency while not deteriorating the separation rate too much.

Similarly to what we did in the previous section, we compare the results obtained with the three different methods on simulated data generated with a parameter $\theta_0 = \sqrt{2 \log(n)/n}$, $\sqrt{\log(n)/n}$ and $0$. The resutls are given in Figure \ref{fig:BFpos}. We observe that the Bayes Factor constructed using the non-local prior of \cite{doi:10.1080/01621459.2015.1130634} has difficulties to detect parameters in $\M_1$ but close to $\M_0$ as positive due to the penalization induced by the prior. On the other hand the usual Bayes Factor based on the simple conjugate Gaussian prior do not detect $\theta=0$ in $\M_0$ while the other two methods have good asymptotic behaviour. 
\begin{figure}[h]
\includegraphics[width = \textwidth]{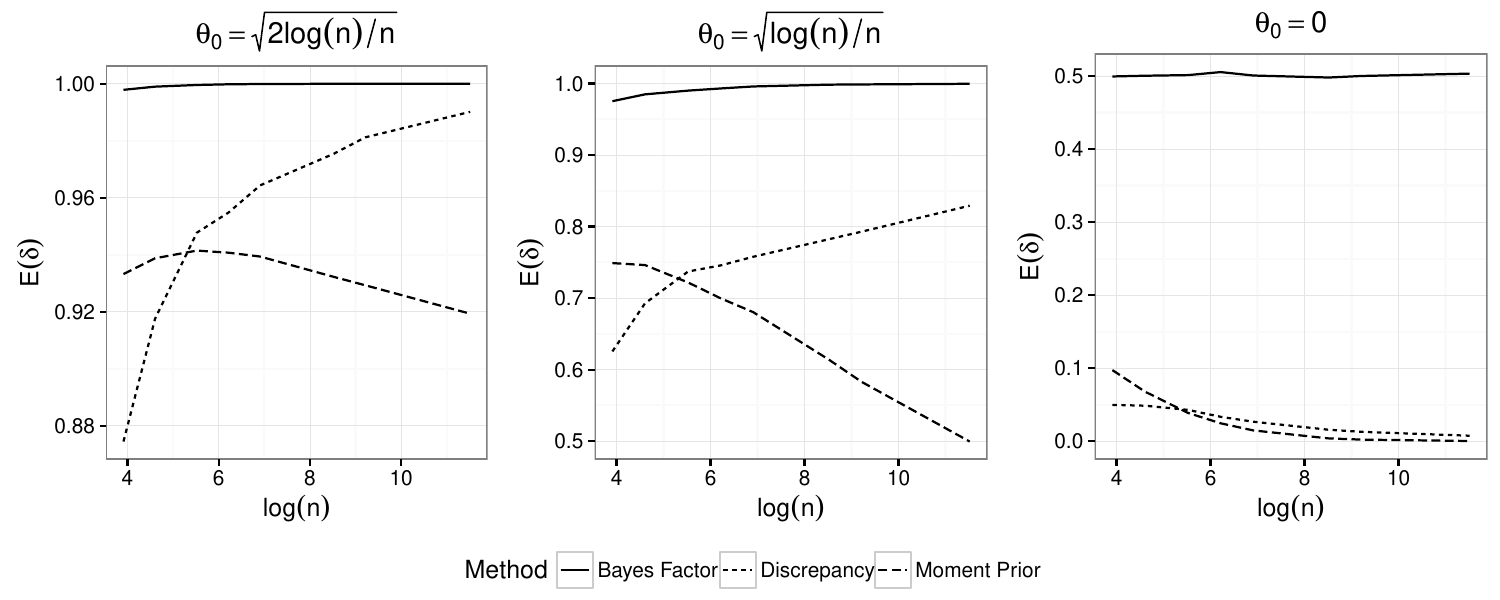}
\begin{center}\caption{Proportion of test that classifies the parameter as negative for $N = 5\times 10^4$ replications of the test. The Bayes Factor obtained with Gaussian and non local priors are compared to 1. For the discrepancy method $\tau_n = \max[1.96,\log(\log(n))]/\sqrt{n}$.}
\label{fig:BFpos}
\end{center}
\end{figure}
We thus see that both the usual Bayes Factor and the approach based on non-local priors have difficulties to detect parameters at or near the boundary. More worryingly, the behaviour of these methods near the boundary strongly depends on the sets $\M_0$ and $\M_1$ and on the prior contraction on these sets. On the other hand the proposed method, although a little less efficient for finite sample sizes, does detect the parameter at or near the boundary. An easy fix in this particular setting to get better results for the Bayes Factor and non local priors would be to elicit $\tilde{\M}_0 = \{0\}$. Nevertheless the same behaviours exposed in the previous section would remain. Furthermore, for more complex hypotheses, it could be difficult to single out the boundary as a separated hypothesis. We shall see in the next section that for these examples, the proposed method attains asymptotically the minimax separation rate. 

\section{Application to standard testing problems}
\label{sec:simple:test}

\subsection{Testing parametric hypotheses}
\label{sec:parametric}

Consider the following parametric model for some fixed $p>0$,
$Y^n \sim P_\theta^n,$ for $\theta \in \Theta \subset \R^p.$
For a fixed subset $\Theta_0 \subset \Theta$ we want to test 
$H_0: \theta \in \M_0 = \Theta_0$, versus  $H_1: \theta \in \M_1 = \Theta \cap \Theta_0^c$. This problem has been widely studied in the Bayesian literature \citep[see][for instance]{Robert2007}. 


In this simple case, the following theorem gives a calibration for the threshold $\tau_n$ in \eqref{eq:test:relaxed} such that the testing procedure satisfies condition \eqref{eq:consistencyH0} and \eqref{eq:consistencyH1}, and gives an upper bound for the separation rate $\rho_n$. 

\begin{theorem}
\label{thm:param}
Let $\Pi$ be a prior distribution on $\Theta$ and $d$ be a metric on the parameter space $\Theta$. Assume that for some positive sequence $\epsilon_n$ we have  
\begin{equation}
\sup_{\theta^* \in \Theta} \E_{\theta^*}^n \Pi[d(\theta,\theta^*) > \epsilon_n | Y^n] = o(1).
\label{eq:theo:param:cvrate}
\end{equation}
Then choosing $D(\theta,\M_{0}) = \inf_{\theta^*\in \Theta_{0}} d(\theta,\theta^*)$ and $\tau_{n} = \epsilon_n$ in \eqref{eq:test:relaxed} the decision rule $\delta^{\pi}_{n}$ in  \eqref{eq:intro:decision:rule} satisfies
\begin{eqnarray*} 
\sup_{\theta \in \Theta_0} \E_\theta^n [\delta^\pi_n(\tau_n)] = o(1), ~
\sup_{\theta \in \Theta, d(\theta,\Theta_0)> 2 \epsilon_n} \E_\theta^n [1-\delta^\pi_n(\tau_n)] = o(1).
\end{eqnarray*}
\end{theorem}

Condition \eqref{eq:theo:param:cvrate} is the standard concentration property of the posterior which is known to hold for regular models with $\epsilon_n = n^{-1/2}u_n$ where $u_n$ is any positive sequence increasing to infinity \citep[see for instance][]{ghosal2000convergence,ghosal:vaart:2007}. In this case the separation rate $\rho_n$ for the proposed test is the minimax separation rate $n^{-1/2}$ up to some factor $u_n$. 
The proof of this Theorem is postponed to Appendix \ref{seq:app:proof:param}. 
From the proof of Theorem \ref{thm:param}, we can also derive an upper bound for $\sup_{\theta \in \Theta_0} \E_\theta^n [\delta^\pi_n(\tau)]$ and $\sup_{\theta \in \Theta, d(\theta,\Theta_0)> 2 \epsilon_n} \E_\theta^n [1-\delta^\pi_n(\tau)]$ of the order of $\sup_{\theta^* \in \Theta} \E_{\theta^*}^n \Pi[d(\theta,\theta^*) > \epsilon_n | Y^n]$. Under some regularity assumptions on the models, we get that the type I and type II error can be uniformly bounded by $e^{-Cu_n^2}$ for some constant $C>0$. Choosing $u_n$ of the order of $\sqrt{\log(n)}$ will thus give polynomial decay uniformly for both errors. As argued in \cite{JohnsonRossell2010}, the Bayes Factor usually contracts at an exponentially fast rate. for a true alternative However, this is to be balanced with the fact that here the proposed control is uniform over all $\theta \in \Theta$ such that $d(\theta,\Theta_0)> 2 \epsilon_n$.

\subsection{Detection of signal in white noise}
\label{sec:WN}
We now apply our approach to the problem of detecting signal in the standard white noise model. 
This problem is closely related to the well studied goodness-of-fit testing problem, where one is interested in testing a parametric hypothesis versus a non parametric one. 
Here again this problem has been extensively studied in the literature. Goodness of fit testing have been considered both from a frequentist and Bayesian point of view, see for instance \cite{ingster2003nonparametric}, \cite{dass2004note} or see \cite{tokdar2010bayesian} for a review. The specific problem of detection of signal in white noise has also been treated in 
 \cite{Ingster87,Lepski1999,Lepski2008}. 

Here we consider the equivalent infinite Gaussian sequence model
\begin{equation}
Y_i = f_i + \frac{\epsilon_i}{\sqrt n}, ~ \epsilon_i \overset{iid}{\sim} \mathcal{N}(0,1), ~ i \geq 1, 
\label{eq:WN:sequence:model}  
\end{equation}
where 
$f = (f_i) \in l^2 = \{g, \sum_i g_i^2 < \infty \}$. 
Similarly to \cite{Lepski1999} we test $f = 0$ against a Sobolev ellipsoid of fixed smoothness $s$, $W_2^s(L) = \{ f \in l_2, \sum_{i=1}^\infty f_i^2 i^{2s} \leq L \}$. We thus have $\M_0 = \{f = 0\}$ and $\M_1 = \{f \in W^2_s(L), f \neq 0\}$. 
 We consider a conjugate Gaussian prior on as in section 3 of \cite{castillo2013general}. For $k_n = n^{2/(4s+1)}$ and all increasing sequence $s = (s_1,s_2,...)$ such that $s_{k_n} \leq n^{4s/(4s+1)}$ and $\sum_{i=1}^{k_n} 1/(n+s_i) \leq \rho_n/4$ we define $\Pi$ by 
\begin{equation}
(f_1, \dots, f_{k_n}) \sim \bigotimes_{i=1}^{k_n} \mathcal{N}(0,s_i^{-1}), ~ f_j = 0 ~ \forall j > k_n. 
\label{eq:WN:prior}
\end{equation}
We choose the discrepancy measure $D(f,\mathcal{M}_{0})$ to be the $l_2$ norm of $f$, $||f||_{2} = (\sum_{i=0}^{\infty} f_{i}^2)^{1/2}$. The following Theorem gives a calibration for the threshold $\tau_{n}$ in \eqref{eq:test:relaxed} and an upper bound for the separation rate of our testing procedure. 
\begin{theorem}
Let $Y^n$ be sample from \eqref{eq:WN:sequence:model} and consider a prior on $f$ as defined in \eqref{eq:WN:prior}. Let $v_n$ be any sequence increasing to infinity and let $\rho_n = v_n n^{-2s/(4s+1)}$ and $\tau_n$ be such that $\tau_n^2 = C\rho_n/2 + k_n/n + \sum_{i=1}^{k_n} \frac{1}{n+s_i}$ for some positive constant $C$. Setting $d$ to be the $l_2$ norm, the decision rule $\delta_{n}^{\pi}$ as defined in \eqref{eq:intro:decision:rule} satisfies 

\begin{equation}
\E_0^n(\delta_n^\pi) = o(1), ~ 
\sup_{f\in \M_1,||f||_2 > \rho_n }  \E_f^n(1-\delta_n^\pi) = o(1) .
\label{eq:result:WN}
\end{equation}

\end{theorem}

Here again the separation rate $\rho_n$ of the test is the minimax separation rate as shown in \cite{Ingster87}. An interesting aspect of this test is that it does not rely on the precise estimation of the true underlying function but rather on the semiparametric estimation of $D(f,\M_0)$ which allows us to obtain a separation rate polynomialy faster than the estimation rate for Sobolev alternative. It is to be noted that the prior \eqref{eq:WN:prior} is not optimal for the estimation problem but leads to the best possible separation rate for the testing problem. The proof of this theorem is postponed to Appendix \ref{sec:app:proof:WN}.  

\section{Shape constraints testing}
\label{sec:shape:constraints}

\subsection{Statistical setting}

We consider the nonparametric fixed design regression problem with Gaussian residuals for $n > 0$
\begin{equation}
Y_j = f(j/n) + \sigma \epsilon_i, ~j=1,\dots,n
\label{eq:shape:model}
\end{equation}
where $\sigma>0$ and $(\epsilon_1,\dots,\epsilon_n)$ is a sequence of independent standard Gaussian random variable. The approach presented in this paper are also valid for non uniform design and random design under additional condition but considering these cases will only make the computations more complex and will thus not be treated here. For this problem, we consider a piecewise constant prior distribution on the regression function $f$ and a prior with density $\pi_\sigma$ with respect to the Lebesgue measure on $\sigma$. More precisely, for $I_{i} = [i-1/k,i/k)$ the uniform partition of $[0,1]$, we define functions $f_{\omega,k}$ as 
$$
f_{\omega,k}(\cdot) = \sum_{i=1}^k \omega_{i}\I_{I_{i}}.
$$ 
We choose the following form for the prior on $f$
\begin{equation}
d\Pi(f) = \pi_k(k) \pi_\omega(\omega_1,\dots,\omega_k|k)d\lambda_k(\omega_1,\dots,\omega_k) d\nu(k),
\label{eq:shape:prior}
\end{equation}
where $\lambda_k$ is the Lebesgue measure on $\R^k$ and $\nu$ the counting measure on $\N$. Note that a similar prior has been studied in \cite{holmes2003generalized} for modelling monotone functions. Here again, although this prior is not well suited for the estimation problem, it gives good theoretical and practical results for testing the shape constraints studied in this paper as shown bellow.
For simplicity we consider a product form for $\pi_\omega$, $\pi_\omega(\omega_1,\dots,\omega_k|k) = \prod_{i=1}^k g(\omega_i)$ where $g$ is a density on $\R$. In addition we assume that the following conditions holds

\begin{description}
     \item[C1] the density $\pi_\sigma$ is bounded and continuous and $\pi_\sigma(\sigma)>0$ for all $\sigma \in (0,\bar{\sigma})$,
     \item[C2] the density $g$ is continuous positive on $\R$ and bounded from above.
     \item[C3] $\pi_k$ is such that there exists positive constants $C_d$ and $C_u$ such
  that
  \begin{equation}
    \label{eq:pik} e^{-C_d kL (k)} \leq \pi_k (k) \leq e^{-C_u kL (k)}
  \end{equation}
  where $L (k)$ is either $\log (k)$ or $1$. 
\end{description}
The condition \textbf{C1} and \textbf{C2} are mild and are satisfied for a large variety of distributions. In section \ref{sec:calibration} we will take $g$ to be a Gaussian density and $\pi_\sigma$ to be a inverse gamma density. Simple algebra shows that for this choice of prior, both conditions are satisfied. Condition \textbf{C3} is a usual condition when considering mixture models with random number of components, see e.g. \citet{Rousseau2010betas} and is satisfied by Poisson or Geometric distribution for instance.

Define the sets 
\begin{align*}
&\FO = \{ f \in L_\infty([0,1]): \forall x \in [0,1], f(x) > 0 \} \\
&\Fm(K) = \{ f: ||f||_\infty \leq K, \forall x \leq y ~f(x) \geq f(y) \}
\end{align*}
of positive and monotone non increasing functions respectively. 
For $\alpha>0$ and $L>0$, define $\mathcal{H}(\alpha,L) = \{f, ||f||_{H,\alpha}\leq L\}$ where $||\cdot||_{H,\alpha}$ is the Hölder norm. 
We consider both testing problems 
$H_0: f \in \M_0 = \FO$ versus $H_1: f \in \M_1 = \mathcal{H}(\alpha,L) \cap \FO^c$ and $H_0: f \in \M_0 = \Fm(K)$ versus $H_1: f \in \M_1 = \mathcal{H}(\alpha,L) \cap \Fm(K)^c$. 
. These problem has been considered in the literature in \cite{Juditsky2002} and \cite{Baraud2005} for instance. 
Note that with a prior chosen as in \eqref{eq:shape:prior} we have $\pi(\FO)>0$ and $\pi(\Fm(K))>0$. Furthermore, if the true regression function $f_0$ is in $\FO$ or $\Fm(K)$ then the piecewise constant function with $k$ pieces of the form \eqref{eq:shape:prior} which minimizes the Kullback Leibler divergence with $P_{f_0}$ will also be in $\FO$, respectively $\Fm(K)$, for all $k$.

We then study the posterior separation rate of the test with respect to the metric $d = d_\infty$ defined as 
\[
d_\infty(f,g) = \sup_{x \in [0,1]} |f(x) - g(x)| .
\]
For each test we compute the the separation rate of our procedure and compare it with the minimax separation rates, which is $n^{-\alpha/(2\alpha +1)}$ in both cases.

Our approach could also apply to other types of shape constraints such as convexity or unimodality using similar methods. 

\subsection{Testing for positivity}
\label{sec:shape:posi}
We first consider positivity constraints. There exist a few methods to test for positivity in a nonparametric setting, see for instance \cite{Baraud2005}. We propose the following discrepancy measure for $D$ in \eqref{eq:test:relaxed} 
\begin{equation} 
D(f,\FO) = -\inf_{x \in [0,1]} f(x).
\label{eq:shape:posi:discrep}
\end{equation}
We immediately have that $D(f,\FO) \leq 0$ if and only if $f \in \FO$. Here the discrepancy measure can be related to the supremum distance with the set of positive functions. For piecewise constant functions $f_{\omega,k}$, $D(f_{\omega,k},\FO)$ has the simple expression $D(f_{\omega,k},\FO) = - \min_{1\leq i \leq k} (\omega_i)$. This turn out to be particularly useful for the calibration of the threshold $\tau_n$. 
Let $\mathcal{G}_k$ be the set of piecewise constant function with $k$ pieces. The idea of the calibration of $\tau_n$ is the following. In the model $\mathcal{G}_k$, the a posteriori uncertainty for estimating $\omega = (\omega_1, \dots,\omega_k)$ is of order $(k/n)^{1/2}$. Hence any positive function $f_{\omega,k}$ such that for all $i$, $\omega_i \geq O\{(k/n)^{1/2}\}$ might be detected as possibly positive. We thus choose a threshold $\tau_n^k$ for each model $\mathcal{G}_k$ of similar order.
The results are presented in the following theorem.

\begin{theorem}
Under the assumptions \textbf{C1} to \textbf{C3}, and if $\underline{\sigma}<\sigma \leq \bar{\sigma}$ for fixed $0<\underline{\sigma} \leq \bar{\sigma}$, then for a fixed constant $M_0 >0$, setting $\tau = \tau_n^k = M_0  \{k \log (n)n^{-1}\}^{1/2}$ and $\delta_n^{\pi}$ the testing procedure defined in \eqref{eq:intro:decision:rule}, for all $K >0$ there exists some $M >0$ such that uniformly for $\alpha \in [\alpha_0,1]$, $\forall \alpha_0>0$
  \begin{equation}
\begin{split}
\sup_{\underline{\sigma}<\sigma \leq \bar{\sigma}} \sup_{f\in \FO}&\E_{f,\sigma}^n(\delta_n^\pi) = o(1) \\ 
\sup_{\underline{\sigma}<\sigma \leq \bar{\sigma}}\sup_{f,d_\infty\{f, \FO\} > \rho , f \in\mathcal{ H}(\alpha,L)}&\E_{f,\sigma}^n(1-\delta_n^\pi) = o(1) 
\end{split}
\label{eq:def_cons2}
\end{equation}
for all $\rho > \rho_n(\alpha) = M\{n/\log(n)\}^{-\alpha/(2\alpha+1)}v_n$ where $v_n = 1$ when $L(k) = \log(k)$ and $v_n = \{\log(n)\}^{1/2}$ when $L(k) = 1$.
    \label{th:main:positivity}
\end{theorem}

\subsection{Testing for monotonicity}

We now consider monotonicity constraints. Tests for monotonicity have been well studied in the frequentist literature, see for instance \cite{Baraud2003,Baraud2005,MR1810919,Bowman1998}. In a Bayesian setting, only \cite{scott2013nonparametric} proposed a test for monotonicity using non local priors. Define the discrepancy measure between $f$ and $\FO$ as 
 
\begin{equation}
D(f,\FO) = \sup_{0\leq x<y\leq 1}\{f(y) - f(x)\}.
\label{eq:shape:mono:discrep}
\end{equation}
Here again when considering piecewise constant functions $f_{\omega,k}$, \eqref{eq:shape:mono:discrep} we get the simple formulation $D(f_{\omega,k},\FO) = \max_{1\leq i \leq j \leq k} (\omega_j - \omega_i) $ which allows for a simple calibration of $\tau_n$ in a similar way as in section \ref{sec:shape:posi}. 

\begin{theorem}
Under the assumptions \textbf{C1} to \textbf{C3}, for a fixed constant $M_0 >0$, setting $\tau = \tau_n^k = M_0  \{k \log (n)n^{-1}\}^{1/2}$ and $\delta_n^{\pi}$ the testing procedure defined in \eqref{eq:intro:decision:rule}, for all $K >0$ 
then there exists some $M >0$ such that uniformly for $\alpha \in [\alpha_0,1]$, $\forall \alpha_0>0$
  \begin{equation}
\begin{split}
\sup_{\underline{\sigma}<\sigma \leq \bar{\sigma}} \sup_{f\in \Fm(K)}&\E_f^n(\delta_n^\pi) = o(1) \\ 
\sup_{\underline{\sigma}<\sigma \leq \bar{\sigma}} \sup_{f,d_\infty\{f, \Fm(K)\} > \rho , f \in\mathcal{ H}(\alpha,L)}&\E_f^n(1-\delta_n^\pi) = o(1) 
\end{split}
\label{eq:def_cons2}
\end{equation}
for all $\rho > \rho_n(\alpha) = M\{n/\log(n)\}^{-\alpha/(2\alpha+1)}v_n$ where $v_n = 1$ when $L(k) = \log(k)$ and $v_n = \{\log(n)\}^{1/2}$ when $L(k) = 1$.
    \label{th:main:monoton}
\end{theorem}
Neither the prior nor the threshold depend on the regularity $\alpha$ of the regression function under the alternative. Moreover for all $\alpha \in (0,1]$, the separation rate $\rho_n(\alpha)$ is the minimax separation rate up to a $\log(n)$ term. Thus our test is almost minimax adaptive. The $\log(n)$ term seems to follow from our definition of the consistency where we do not fix a level for the Type I or Type II error contrariwise to the frequentist procedures. The conditions on the prior are quite loose, and are satisfied in a wide variety of cases. The constant $M_0$ does not influence the asymptotic behaviour of our test but has a great influence in practice for finite $n$. A practical way of choosing $M_0$ is given in section \ref{sec:calibration}.

\section{Simulation study for positivity and monotonicity testing}

\subsection{Prior specification and sampling strategy}
\label{sec:calibration}

Conditions on the prior in Theorem \ref{th:main:monoton} are satisfied for a wide variety of distributions. However, when no further information is available, some specific choices can ease the computations and lead to good results in practice. We present in this section such a specific choice for the prior and a way to calibrate the hyperparameters. We also fix $\gamma_0 = \gamma_1 = 1/2$ in the definition of $\delta_n^\pi$.

A practical default choice is the usual conjugate prior, given $k$, i.e. a Gaussian prior on $\omega$ with variance proportional to $\sigma^2$ and an Inverse Gamma prior on $\sigma^2$. This will considerably accelerate the computations as sampling under the posterior is then straightforward. Condition \eqref{eq:pik} on $\pi_k$ is satisfied by the two classical distributions on the number of parameters in a mixture model, namely the Poisson distribution and the Geometric distribution. It seems that choosing a Geometric distribution is more appropriate as it is less spiked. We thus choose for $\lambda,a,b>0$, $m\in \R$ and $\mu >0$

\begin{equation}
    \Pi= 
    \begin{cases}
        k \sim \mathrm{Geom}(\lambda) \\
        \sigma^2 | k \sim IG(a,b) \\
        \omega_i | k, \sigma \overset{iid}{\sim} \mathcal{N}(m,\sigma^2/\mu)  
    \end{cases}
    \label{eq:defpi}
\end{equation}
Standard algebra leads to a close form for the posterior distribution up to a normalizing constant. Let $n_i =\mathrm{Card} \left\{j, j/n  \in [(i-1)/k,i/k)  \right\}$, we denote 
\[ 
\tilde{b}_k = b + \frac{1}{2}  \sum_{i = 1}^k \left\{ \sum_{j, j /
   n \in I_i} \left( Y_j - \overline{Y_i} \right)^2 + \frac{n_i \mu}{n_i +
   \mu} ( \overline{Y_i} - m)^2 \right\}, 
\]
where $\overline{Y_i}$ is the empirical mean of the $Y_l$ on the set $\{j, j /
n \in [(i-1)/k,i/k) \}$, we have
\[ \pi_k (k|Y^n) \propto \pi (k)  \tilde{b}_k^{- (\alpha + n / 2)} \mu^{k
   / 2}  \prod_{i = 1}^k (n_i + \mu)^{- 1 / 2}
\]
We can thus compute the posterior distribution of $k$ up to a constant. 
We will thus be able to sample from $\pi_k(k|Y^n)$ using a truncated approximation of the posterior. In the examples we choose to truncate at some $k_0\leq n$. 
We then compute the posterior distribution of $\omega$ and $\sigma$ given $k$
\begin{align*}
\sigma^2 | k,Y^n &\sim IG(a + n/2,\tilde{b}_k) \\
\omega_j | k,\sigma^2, Y^n &\overset{ind.}\sim \mathcal{N}\left({m\mu + n_j \bar{Y}_j \over n_j + \mu}, {\sigma^2 \over n_j + \mu} \right)  .
\end{align*} 
Given $k$, sampling from the posterior is thus straightforward. 

A crucial hyperparameter that needs to be calibrate is for $M_0$ the constant in $\tau$. A close inspection of the proofs (in particular the proof of Lemma \ref{lem:pi1/2}) using the fact that we have a Gaussian posterior, gives us that taking 
$$\tau_n = \sqrt{\frac{\log(k/n) k \sigma^2}{n+k\mu \sigma^2}},$$ would induce the desired results. 

\subsection{Simulated Examples}
In this section we run our testing procedure on simulated data to study the behaviour of our test for finite sample sizes. We first examine the behaviour of the proposed test for positivity on an example that illustrate that the separation rate of the test is indeed upper bounded by $(\log(n)/n)^{\alpha/(2\alpha +1)}$ up to some constant. 
We then compare our test for monotonicity to other methods proposed in the literature, and get comparable results for finite sample size.  
\subsubsection{Testing for positivity}

Consider the test for positivity proposed in section \ref{sec:shape:posi}. Similarly to the examples of section \ref{sec:boundary}, we will consider a sequence of function that are in $\M_1$ i.e. not positive, but are getting closer and closer to the boundary. More precisely we take 
$$
f_n(x) = 10\rho_n(|x-0.1| - 0.1)\I_{|x-0.5|<0.1},
$$
and thus $\rho_n = d_\infty(f,\F_{+})$. Plots of $f_n$ for different values of $n$ are given in Figure \ref{fig:funcpos}.
\begin{figure}[h]
\begin{center}
\includegraphics[scale = 0.5]{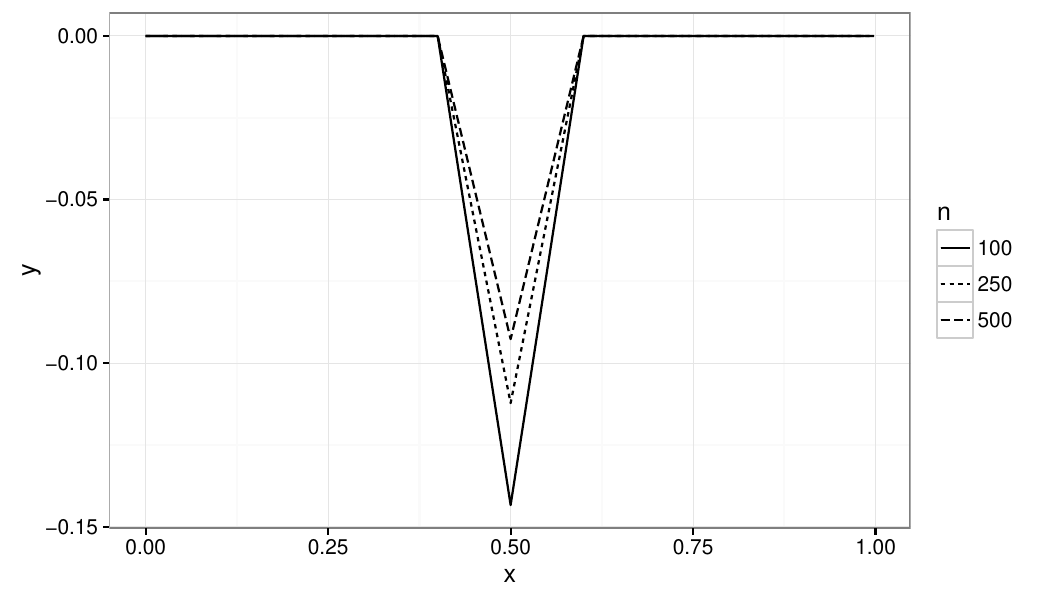}
\caption{\label{fig:posres} Plots of $f_n$ for different values of $n$ with $\rho_n = 0.4 (\log(n)/n)^{1/3}$}
\end{center}
\end{figure}
Since for all $n$ this function is piecewise linear, we thus have that $f \in \mathcal{H}(\alpha,L)$, with $\alpha = 1$. Given Theorem \ref{th:main:positivity}, we have that for some constant $M$ large enough, the test should be consistent for $f_n$ if $\rho_n > M(\log(n)/n)^{1/3}$. 

We run our test on simulated data generated from the model \eqref{eq:shape:model} with $f=f_n$ for different values of $M$ and with $f=0$ that lies at the boundary between hypotheses. The results are given in figure \ref{fig:posres}. We observe that the test detects parameter at the boundary as positive, even for moderate values of $n$. In addition, for $M>0.4$, the function $f_n$ are detected as non-positive, and the asymptotic regime is attained around $n = 2000$, while for $M< 0.4$ the functions $f_n$ are not detected as non-positive. This indicates that the test does separate the hypotheses at the rate at least $0.4(\log(n)/n)^{1/3}$, and we thus recover the results from Theorem \ref{th:main:positivity}.    
\begin{figure}[h]
\begin{center}
\includegraphics[scale = 0.55]{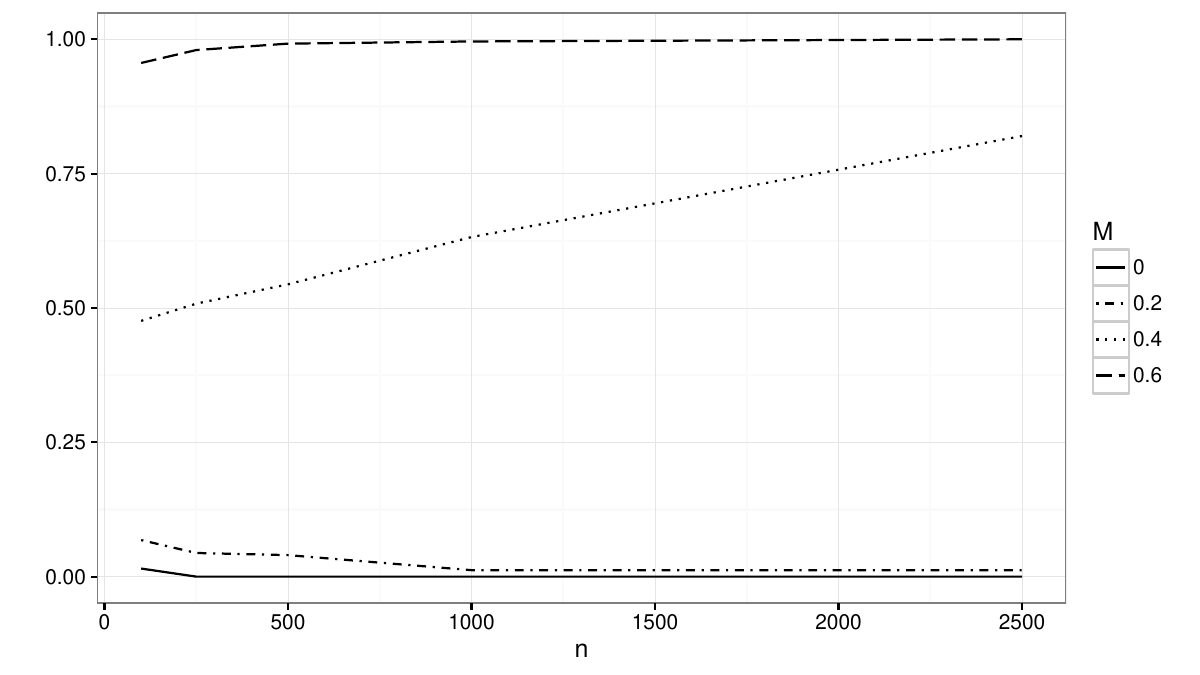}
\caption{\label{fig:funcpos} Proportion of functions classified as non-positive estimated on $K=500$ independent replication of simulated data generated from model \eqref{eq:shape:model} with $f = 10 M (\log(n)/n)^{1/3}(|x-0.1| - 0.1)\I_{|x-0.5|<0.1}$ for different values of $M$.}
\end{center}
\end{figure}

\subsubsection{Testing for monotonicity}
\label{sec:simu}
We now compare our approach to test for monotonicity with the ones proposed in the literature. We consider the following nine functions adapted from \citet{scott2013nonparametric} and \citet{Baraud2003} and plot in Figure \ref{fig:func}. 

\begin{equation} 
\begin{aligned}[l]
    f_1 (x) = & - 4 (x - 0.5)^3 \I_{x \leq 1 / 2}  - \\& 0.1 (x - 0.5) + 0.25e^{- 250
    (x - 0.25)^2}\\
    f_2 (x) = & 0.1 x\\
    f_3 (x) = & 0.1 e^{- 50 (x - 0.5)^2}\\
    f_4 (x) = & - 0.1 \cos (6 \pi x)\\
    f_5 (x) = & - 0.2 x + f_3 (x)\\
    \end{aligned} 
    \begin{aligned}[l]
    f_6 (x) = & - 0.2 x + f_4 (x)\\
    f_7 (x) = & - (1 + x) + 0.25 e^{- 50 (x - 0.5)^2}\\
    f_{8}(x) = & -x -1 + 0.45e^{-50(x-0.5)^2} \\
    f_9 (x) = & - 0.5 x^2\\
    f_{10} (x) = & 0\\
    f_{11}(x) = & -x-1 
    \end{aligned}
  \label{eq:functions}
\end{equation}
\begin{figure}[h]
  \includegraphics[width = \textwidth]{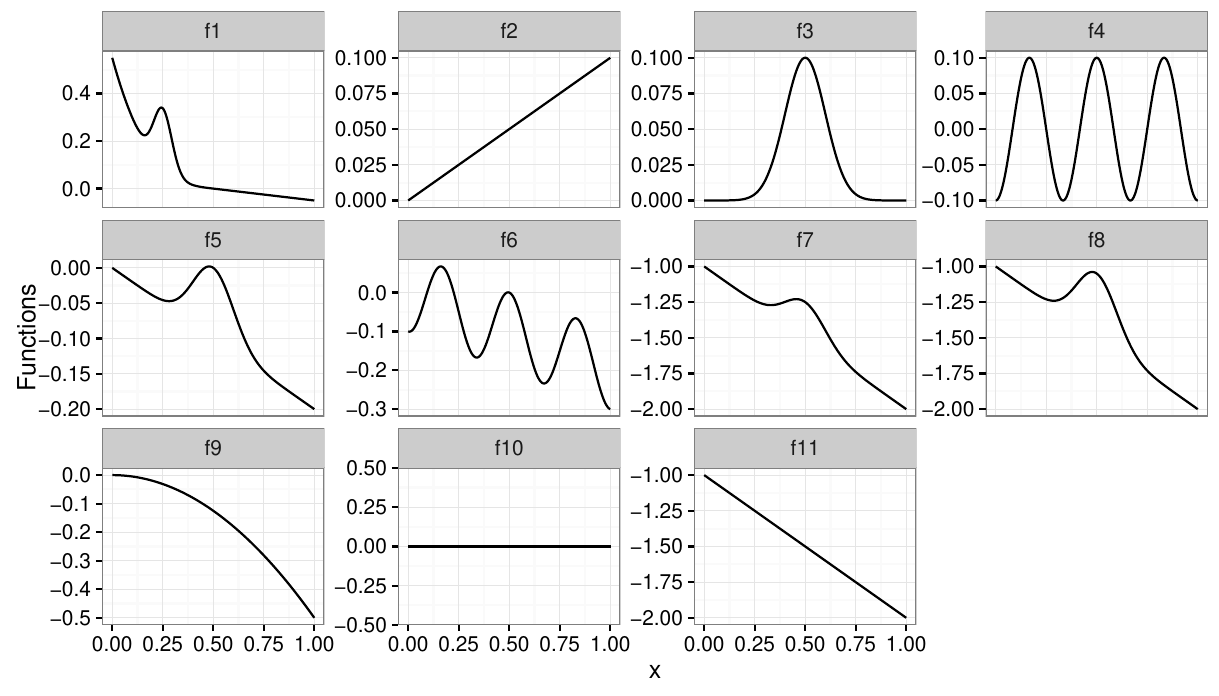}
    \caption{Regression functions used in the simulated examples.}
\label{fig:func}
\end{figure}
The functions $f_1$ to $f_6$ are clearly not in $\Fm(K)$ with $K=2$. The function $f_7$ has a small bump around $x = 0.5$ which can be seen as a local departure from monotonicity. This function is thus expected to be difficult to detect for small datasets given our parametrization. The function $f_9$ is a completely flat function and belongs to $\Fm(K)$.    

For several values of $n$, we generate $N = 500$ replicates of the data $Y^n = \{y_i, i=1,\dots ,n\}$ from model \eqref{eq:shape:model}. For each dataset, we approximate ${\pi}\left\{ D(f_{\omega,k},\Fm) > \tau_n^k | Y^n \right\} $ based on $K = 5\times 10^4$ samples from the posterior and reject the null if 
\[
\hat{\pi}\left\{ D(f_{\omega,k},\Fm) > \tau_n^k | Y^n \right\} = {1 \over K} \sum_{i=1}^K \I\left\{ D(f_{\omega^i,k^i},\Fm) > \tau_n^{k^i} \right\} \geq 1/2
\]
The results are given in table \ref{tab:res:simu}. \\ 

\begin{table}[h]
{\begin{center}

\begin{tabular}{c|c|c||ccc|cc}
&{$f_0$}& {$\sigma^2$}&{Barraud}& {Akakpo} & {Scott} & {Discrepancy}  \\ 
  \hline
  \multirow{8}{2ex}{$H_1$}  &    
   $f_1$ &0.01  &6   &9  &100 &80   \\
&  $f_2$ &0.01  &64  &33 &74  &75  \\
&  $f_3$ &0.01  &53  &43 &35  &77  \\
&  $f_4$ &0.01  &92  &92 &91  &100  \\
&  $f_5$ &0.01 &24  &25 &85  &17  \\
&  $f_6$ &0.01 &77  &75 &99  &99  \\
&  $f_7$ &0.01  &1   &4  &91  &39  \\
&  $f_8$ &0.01  &71  &82 &99  &99   \\
\hline
\multirow{3}{2ex}{$H_0$} 
             &  $f_9$    &0.01  & 100 &100 & 93 &93 \\
             &  $f_{10}$ &0.01  & 97  &94  & 95 &89  \\
             &  $f_{11}$ &0.01  & 100 &99  & 97 &94  \\ 
\hline
\multicolumn{3}{c||}{Average} &62.3 & 59.5 & 87.2 & 78.36  
\end{tabular}
\caption{\label{tab:res:simu} Results of the simulation study. Each entry is the percentage of correctly classified functions estimated on $K = 500$ independent replication of the experiment. Baraud method for \cite{Baraud2003}, Akakpo for \cite{ABD2012}, Scott for \cite{scott2013nonparametric} with Gaussian prior and Discrepancy  for the proposed method.}
\end{center}
}
\end{table}
For all the considered functions, the computational time is reasonable even
for large values of $n$. For instance, for $f_1$, we require less than $2$
seconds to perform the test for $n = 2500$ using a simple R script available on demand. We compare our results with the ones obtained in \citet{scott2013nonparametric} for the Gaussian prior and the methods proposed by \cite{Baraud2003} and \cite{ABD2012}. The results are given in Table \ref{tab:res:simu}. The proposed method is a little less efficient than the one based on non-local prior in average, but it seems to perform better for some functions (e.g. $f_3$). When $n$ grows, the percentage of correctly classified function goes to 1 as predicted by the theory.

\section{Proofs}

\subsection{Proof for the parametric test}
\label{seq:app:proof:param}

We prove that with the proposed calibration for $\tau_n$ and $D(\theta,\Theta_0) = \inf_{t \in \Theta_0} d(\theta,t)$, the decision rule satisfies \eqref{eq:consistencyH0}-\eqref{eq:consistencyH1}. For $\theta^* \in \Theta_0$, we directly have that
$$
D(\theta,\Theta_0) \leq d(\theta,\theta^*)
$$
which together with \eqref{eq:theo:param:cvrate} gives \eqref{eq:consistencyH0}. Now for $\theta^* \in \Theta \cap \Theta_0^c$ we have if $d(\theta^*, \Theta_0) > \rho_n$
$$
D(\theta,\Theta_0) \geq d(\theta^*,\Theta_0) - d(\theta,\theta^*) \geq \rho_n - d(\theta,\theta^*).
$$
We deduce \eqref{eq:consistencyH1} directly from condition \eqref{eq:theo:param:cvrate} which ends the proof. 
\subsection{Proof for the detection of signal in white noise}
\label{sec:app:proof:WN}

We fist prove that with the proposed calibration of $\tau_n$ the decision rule \eqref{eq:intro:decision:rule} satisfies \eqref{eq:consistencyH0}. In the sequel $c$ will denote a generic absolute constant that may change from one line to another.   
We want to bound $\Pi(||f||^2 > \tau_n^2 | Y^n)$ when $f_0 = 0$ for $\tau_n^2 = \rho_n/2 + k_n/n + \sum_{i=1}^{k_n} \frac{1}{n+s_i} $.
For all $t \leq 2 n$ we have, using the Chernoff bound we have with $P_0^n$-probability that goes to $1$ 

\begin{align*}
\log\{ \Pi(||f||^2 > \tau_n^2 | Y^n)\}  &\leq  -t \tau^2_n + \sum_{i=1}^{k_n} t \frac{\epsilon_i^2}{n} + \sum_{i=1}^{k_n} \frac{1}{n+s_i} \\
&\leq -t\frac{\rho_n}{2} \\ 
&\leq 1-c  
\end{align*}
for $c$ large enough, which give the result

We now state an auxiliary result that will be needed for the remainder of the proof. Define $H(s,\rho) = \{f \in W_2^s(L), ||f||_2 > \rho \}$

\begin{lemma}
Let $k_n = n^{2/(4s +1)}$ and consider $Z_n = \sum_{i=1}^{k_n} (Y_i^2 - 1/n)$. We thus have if $f_0 \in H(s, \rho_n)$ with $\rho_n =v_n n^{-4s/(4s+1)}$ where $v_n \to \infty$ slowly with $n$ then for some $C>0$

\begin{equation}
P_0^n(Z_n \leq \rho_n ) = o(1)
\end{equation}
\end{lemma}
The proof of this Lemma can be found in the supplementary materials. 

We now end the proof by showing that $\delta_n^\pi$ satisfies \eqref{eq:consistencyH1}. 
We want to bound $\Pi(||f||^2 \leq \tau_n | Y^n)$ when $f_0 \in H(s, \rho_n)$. For all $n/2>t>0$  and all increasing sequence $s_i$ such that $s_{k_n} \leq n^{4s/(4s+1)}$ we have for $Y^n$ such that $P_0^n(Z_n > \rho_n )$, using the Chernoff bound and the fact that $\sum_{i=1}^{k_n}1/(n+s_i) \leq \rho_n/4$

\begin{align*}
\log\{\Pi(||f||^2 \leq \tau_n^2 | Y^n)\} & \leq  tu_n - \sum_{i=1}^{k_n} \left\{ \frac{tn^2 Y_i^2}{(n+s_i) (n+s_i + 2t)} - \frac{1}{2} \log\Big(1 + \frac{2t}{n+s_i}\Big) \right\} \\ 
&\leq  t\tau_n^2 -  \sum_{i=1}^{k_n} \left\{ tY_i^2  (1 - 2\frac{s_i + t}{n}) - \frac{t}{n+s_i} \right\} \\ 
& \leq  - \frac{t\rho_n}{4} +  2 k_n \frac{s_{k_n} + t}{ n^2} + \frac{k_n t^2}{n^2}  ~  \\
& \leq c ,
\end{align*}
for some $c$ for $v_n$ large enough by taking $t \asymp \rho_n^{-1} $ and $s_{k_n} \leq t^2$,
where the second line comes from the fact that $\frac{1}{(n+s_i)(n+s_i+2t)} \geq \frac{1}{n^2}\Big( 1 - 2\frac{s_i + t}{n} \Big)$.

\subsection{Proof for shape constraints}
\label{sec:app:proof:shape}
\subsection{Auxiliary result}
For all functions $f_0$ in $L_\infty([0,1])$ denote by $P_0$ the probability distribution of $Y^n$ generated with $f = f_0$ and $f_{\omega^0,k}$ the function of $\mathcal{G}_k$ the set of piecewise constant functions with $k$ pieces, that minimizes the Kullback Leibler divergence between $P_f$ and $P_0$. Standard computation gives 
\begin{equation}
\omega_i^0 = n_i^{-1} \sum_{j, j/n \in [(i-1)/k,i/k)} f_0(j/n) , \; n_i = \mathrm{Card}\left\{ j, j/n \in [(i-1)/k,i/k) \right\}. 
\label{eq:def:omega0}
\end{equation} 
The following lemma gives some concentration result for $f_{\omega,k}$ that will be useful for the study of $D(f,\FO)$ or $D(f,\Fm(K))$ respectively for both monotonicity and positivity constraints. 
 \begin{lemma}
     Let $M$ be a positive constant. Let $\Pi$ be as define in \eqref{eq:shape:prior} such that it satisfies condition \textbf{C1}, \textbf{C2} and \textbf{C3}. Denote by $\omega_0$ the minimizer of the Kulback-Leibler divergence $KL(P_{f_{\omega,k}}, P_0)$. Then if there exists a constant $C$ such that $\Pi(\sigma_0/\sigma < C|Y^n) = o_{P_0^n}(1)$ for a constant $A>0$ large enough, we have  
\begin{equation} 
P_0^n \left\{ \Pi\left( \max_{j=1,\dots,k}|\omega_j - \omega^0_j| \geq A \xi_n^k |Y^n \right) \leq \frac{\gamma_1}{\gamma_0 + \gamma_1} \right\} \to 1. 
\label{eq:pi:1/2} 
\end{equation}
where $\xi_n^k = \left[\{k\log(n)\}/{n}\right]^{1/2}$ for all fixed positive $\gamma_0$ and $\gamma_1$.
\label{lem:pi1/2} 
 \end{lemma}
 The proof of this lemma is given in the supplementary materials.
We also state the following lemma that gives a control on the posterior distribution of $k$. 
\begin{lemma} 
Let $k_n = n\epsilon_n^2/\log(n)$ if $L(k) = \log(k)$ and $k_n = n \epsilon_n^2$ if $L(k) = 1$ where $\epsilon_n$ is either $\epsilon_n(\mathcal{F})$ if $f_0 \in \F$ or $\epsilon_n(\alpha)$ if $f_0 \in \mathcal{H}(\alpha,L)$. For $C_1$ a positive constant that my depend on $K$ or $L$, let $\mathcal{K}_n = \{ k \leq C_1 k_n \}$. If $\Pi$ is define as in \eqref{eq:shape:prior} and satisfies \textbf{C1} or \textbf{C1'}, \textbf{C2} and \textbf{C3} we have 
\begin{equation}
\Pi\left( \mathcal{K}_n^c|Y^n \right) \leq o_{P_0^n}(1) 
\label{eq:pikn}
\end{equation}
\label{lem:pikn}
\end{lemma}
The proof is given in the supplementary materials.

\subsubsection{Proof for the test for positivity}

We first prove that $\delta_n^\pi(\tau_n)$ satisfies \eqref{eq:consistencyH0} for $\tau_n = \{k\log(n)/n\}^{1/2}$. Let $f_0 \in \FO$, then for all $k>0$ we have $f_{\omega^0,k} \in \FO$ which in turns gives

$$
D(f_{\omega,k},\FO) = -\min(\omega) \leq \max_{j=1,\dots,k} |\omega_j - \omega^0_j|. 
$$

Note that if $\sigma_0 \leq \bar{\sigma}$ we get directly for $C$ large enough that $\Pi(\sigma_0/\sigma < C | Y^n) = o_{P_0^n}(1)$

Applying lemma \ref{lem:pi1/2} gives us immediately \eqref{eq:consistencyH0} for $M_0$ large enough. We now show that $\delta_n^\pi(\tau_n)$ satisfies \eqref{eq:consistencyH1} with $\rho = \rho_n = M \{ n/\log(n) \}^{-\alpha/(2\alpha+1)} v_n$ for $v_n$ as in Theorem \ref{th:main:positivity}. First note that for $f_0$ such that $f_0 \in \mathcal{H}(\alpha,L)$ $d_\infty(f,\FO) > \rho_n$ we have for all $k$
$$
-\min_{j=1,\dots,k}(\omega^0_j) \geq \rho_n - k^{-\alpha}, 
$$

which leads to 
\begin{align*}
- \min_{j=1,\dots,k} (\omega_j) &\geq -\min_{j=1,\dots,k}(\omega_j^0) - \max_{j=1,\dots,k}|\omega_j - \omega_j^0| \\ 
& \geq \rho_n - k^{-\alpha} - \max_{j=1,\dots,k}|\omega_j - \omega_j^0|. 
\end{align*}
We thus deduce the following upper bound for $\Pi\{D(f,\FO) \leq \tau_n|Y^n\}$: 
$$
\Pi\{D(f,\FO) \leq \tau_n|Y^n\} \leq \Pi(\max_{j=1,\dots,k}|\omega_j - \omega_j^0| \geq \rho_n - k^{-\alpha} - \tau_n | Y^n). 
$$
We ends the proof by applying Lemma \ref{lem:pi1/2} together with Lemma \ref{lem:pikn}.

\subsubsection{Proof for the test for monotonicity}

We first prove consistency under $H_0$. Let $f_0 \in \F$ then
  \[ D(f_{\omega,k},\Fm) \leq 2 \max_{i=1,\dots,k} | \omega_i - \omega_i^0 | \]
  and thus
  \[ P_0^n \left[ \Pi \{D(f_{\omega,k},\Fm) \geq  \tau_n^k |Y_n\} < \frac{\gamma_1}{\gamma_0 + \gamma_1} \right]\to 1 \]
 as soon as $\tau_n^k \geq 2 A \xi_n^k$, which gives the consistency under $H_0$ given Lemma \ref{lem:pi1/2}.

We now prove consistency under $H_1$.
  Let $f_0 \not\in \F$ and $f_0 \in \mathcal{H}(\alpha, L)$ we have
  \begin{equation}
    D(f_{\omega,k},\Fm) \geq D(f_{\omega^0,k},\Fm) - 2 \max_{i=1,\dots,k} | \omega_i - \omega_i^0 |
    \label{eq:hom0}
  \end{equation}
  Assume that $\rho_n(\alpha) < d_\infty (f_0, \F)$, we derive a lower bound for $D(f_{\omega^0,k},\Fm)$. Let $g^{\ast}$ be the monotone non increasing piecewise constant function on the partition
  $\{[0,1/k), \dots, [(k-1)/k,1)\}$, with for $1 \leq i \leq k$, $g^{\ast}_i = \min_{j \leq i}
  \omega^0_j$. Given that $d_\infty(f_{\omega^0,k}, \F) = \inf_{g\in \F} d_\infty(f_{\omega^0,k}, g)$ we get 
  \[ d_\infty (f_{\omega^0, k}, \Fm) \leq d_\infty( f_{\omega^0, k}, g^{\ast})  \leq D(f_{\omega^0,k},\Fm) \]
     And therefore, given that $d_\infty(f_0,\Fm) \leq d_\infty(f_{\omega^0,k},\Fm) + d_\infty(f_{\omega^0,k}, f_0)$ 
  \[ \Pi \left\{ D(f_{\omega,k},\Fm) < \tau_n^k |Y_n \right\} \leq \Pi \left\{ \max_{i=1,\dots,k} |
     \omega_i - \omega_i^0 | \geq \frac{\rho_n(\alpha) - d_\infty(f_{\omega^0, k} , f_0) 
     - C\tau_n^k}{4} |Y^n \right\} \]
  
  For $\mathcal{K}_n$ as in Lemma \ref{lem:pikn} and $k \in \mathcal{K}_n$ and $M$ large enough we have $\rho_n(\alpha) / 4 >
  \tau_n^k$. 
  On the set $\mathcal{K}_n$ we have for $M$, the constant in $\rho_n(\alpha)$ large enough $\rho_n(\alpha) /
  4 \geq d_\infty( f_{\omega^0, k},f_0)$ which in turns gives 

  \[ \Pi \left\{ D(f_{\omega,k},\Fm) < \tau_n^k |Y_n \right\} \leq \Pi \left[ \{\max_{i=1,\dots,k} |
     \omega_i - \omega_i^0 | \geq \rho_n(\alpha) / 8\} \cap \{\mathcal{K}_n
     \cap B_n \}| Y^n \right]  + o_{P_0^n} (1) . \]
  Given \eqref{eq:pi:1/2}, we get that for all $f_0$ such that $d_\infty (f_0, \Fm)
  > \rho_n(\alpha)$
  \[ P_0^n\left[ \Pi \{ D(f_{\omega,k},\Fm) < \tau_n^k |Y_n\} < \frac{\gamma_0}{\gamma_0 + \gamma_1} \right] \to  1 \]
  which ends the proof.

\section{Discussion}

In this paper we present an approach for testing un-separated hypotheses that relies on the estimation of a distance between the parameter and the null set. This approach can be viewed either as a modification of the testing loss function or as a relaxation of the hypotheses at hand. The test obtained using this approach have been shown to be consistent and to achieve the minimax separation rates when testing parametric hypotheses and in some nonparametric settings. 

The approach proposed here currently focus on two hypotheses testing, however, we believe that it could also be applied to more general settings, in particular for the sparse Gaussian sequence model. In this case \cite{MR2650751} proposed a model selection method for the Horseshoe prior. The idea of their approach is somehow similar to the one proposed here. 
\cite{MR2850212} and \cite{MR3036256} have derived some upper bounds on the multiple testing risk from asymptotic properties of each individual test for the Horseshoe prior. We thus believe that the approach presented in this paper could also lead to interesting results in this setting.
In particular one could adopt a minimax version of the risk studied in \cite{MR2850212} and adapt the approach studied here.

\bibliographystyle{ba}
\bibliography{test_final_mref}

\begin{acknowledgement}

I am grateful to the Editor, Associate Editor and the two anonymous referees for careful reading and suggestions. I also thank Judith Rousseau and Peter Gr\"unwald for helpful discussion.

\end{acknowledgement}

\appendix

\section{Proof of Lemma 1}

First recall that $Z_n = \sum_{i=1}^{k_n} (Y_i^2 - 1/n)$ and note that 
\begin{multline*}
P_0^n(Z_n \leq C \rho_n ) =  P_0^n \left\{ ||f_0||^2 - \sum_{i=k_n+1}^\infty f_{0,i}^2 + \frac{2}{\sqrt n} \sum_{i=1}^{k_n} \epsilon_i f_{0,i} + \frac{1}{n} \sum_{i=1}^{k_n} (\epsilon_i^2 - 1)  \leq C \rho_n \right\} \\ 
\end{multline*}
Furthermore, we have that 
\begin{align*}
\E \left\{ \left( \frac{2}{||f_0||^2 \sqrt{n}} \sum_{i=1}^{k_n} \epsilon_i f_{0,i} \right)^2 \right\} &\leq \frac{4}{n||f_0||^2} \leq \frac{4}{n\rho_n} = o(1).
\end{align*}
We thus have that $\frac{2}{\sqrt n} \sum_{i=1}^{k_n} \epsilon_i f_{0,i} = o_{P_0}( ||f_0||^2 )$. We deduce that for $C_0$ large enough

$$
P_0^n \left( Z_n \leq C \rho_n \right) = \P_0^n \left\{ \frac{1}{n} \sum_{i=1}^{k_n} (\epsilon_i^2 - 1) \leq - C' \rho_n \right\} = o(1).
$$
given that $ k_n^{-1/2} \sum_{i=1}^{k_n} (\epsilon_i^2 - 1)$ is asymptotically standard Gaussian. 

\section{Existence of test for the regression model with unbounded variance for the residuals}

We prove the existence of an exponentially consistent sequence of test when estimating the mean of a Gaussian vector with unbounded variance. This result has an interest in its own as it extend the existing results on Bayesian nonparametric regression developed in \cite{ghosal:vaart:2007} or \cite{deJongeZantenBVM}.

Consider the nonparametric regression problem 
$$
Y_i = f_i + \sigma \epsilon_i,~\epsilon_i \overset{iid}{\sim} \mathcal{N}(0,1)  ~i=1,\dots, n 
$$

Let $\Pi$ be a prior on $f$ and let $d_n(f,g) = n^{-1} \sum (f_i - g_i)^2$. 
We have the following lemma 

\begin{lemma}
For $\F_n$ a sequence of sieves, define $\F_n^j = \{f \in \F_n, j^2 \epsilon_n^2 \leq d_n(f,f_0)^2 + (\sigma - \sigma_0)^2 \leq (j+1)^2 \epsilon_n^2 \}$
and assume that for all $j$ we can have a $\epsilon_n$-net for $\F_n^j$ with at most $e^{Cjn\epsilon_n^2/2} $ points. Then there exists a sequence of tests $\Phi_n$ such that 
$$
P^n_0 \Phi_n \to 0, ~ \sup_{f \in \F_n^j} P^n_f(1-\Phi_n) \leq e^{-Cj^2 n \epsilon_n^2/2}
$$
\label{lem:tests}
\end{lemma}

\begin{proof}

We consider 3 cases : $|\sigma - \sigma_0| \leq 1/2$, $\sigma > 3\sigma_0/2$ and $\sigma\leq \sigma_0/2$. \\

When $|\sigma - \sigma_0| \leq 1/2$: we can construct a test $\Psi_1$ such that

\[
\E_0^n(\Psi_1) \leq e^{-Cj^2 n \epsilon_n^2} ; \sup_{ \F_j^k \cap \left\{ |\sigma - \sigma_0| \leq \sigma_0/2 \right\} } \E_{f,\sigma} (1- \Psi_1) \leq e^{-Cj^2 n \epsilon_n^2}. 
\]
as in Lemma 2 of \cite{ghosal:vaart:2007} since the $d_n$ norm can be related to the Hellinger distance in this case. 

For $\sigma > 3\sigma_0/2$ we consider the test $\Psi_2$ defined as
\[
\Psi_2 = \I\left\{ \sum_{i=1}^n \left( \frac{Y_i - f_{0,i}}{\sigma_0} \right)^2 > n c_1 \right\} ,
\]
for a suitably choosen constant $c_1>0$. Chernoff bound gives
\[
\E_0^n (\Psi_2) \leq e^{-Cn} .
\]
If $\sigma> 3 \sigma_0/2$ and $(f,\sigma) \in \F_n^j$, thus $j > j_0/\epsilon_n$ for some $j_0>0$. If $Y_i = f_i + \sigma \eps_i$ where $\eps_i \sim \mathcal{N}(0,1)$ then $ \sum_{i=1}^n  \left( {(Y_i - f_{0,i})}/{\sigma_0} \right)^2 $ follow a non central $\chi^2_n$ distribution with non centrality parameter $\sum_{i=1}^n (f(x_i) - f_{0,i})^2 / \sigma^2 >0$. Thus setting $W \sim \chi^2_n$  
\[ 
\E_{f,\sigma}(1- \Psi_2) = P_{f,\sigma}\left[ \frac{\sigma^2}{\sigma_0^2} \sum_{i=1}^n  \left\{ \frac{Y_i - f_{0,i}}{\sigma} \right\}^2 \leq n c_1 
\right] \leq \mathrm{pr}\left( W \leq \frac{4}{9} c_1n \frac{\sigma_0^2}{\sigma} \right).
 \]
Chernoff bound gives 
\[
E_{f,\sigma}(1 - \Psi_2) \leq e^{-C_2 n} .
\]

Recall that we can construct a $\epsilon$-net for $\F_j$ with less that $e^{Cj^2\epsilon^2/2}$ points.  
For $\sigma < \sigma_0/2$ we consider the test $\Psi_3^t$ associated to $f^t \in \F_n^j$ a point in the $\xi \epsilon_n$ net and some suitably chosen $0< c_2 < 1 $ defined as 
 \[
 \Psi_3^t = \I\left[ \sum_{i=1}^n \left\{ \frac{Y_i - f^t_i}{\sigma_0} \right\}^2 \leq c_2 n \right].
 \]

 As before, given that under $P_{f_0, \sigma_0}$, $\sum_{i=1}^n \left[ \{Y_i - f^t_i\}/{\sigma_0} \right]^2$ follows a non central $\chi^2_n$ distribution

\[
\E_0^n (\Psi_3^t) = P_0 \left[ \sum_{i=1}^n \left\{ \frac{Y_i - f^t_i}{\sigma_0}\right\}^2 \leq c_2 n \right] \leq \mathrm{pr}(W \leq c_2n).
\]

Given that the moment generating function of a non central $\chi^2_n$ distribution with non centrality parameter $\Delta$ at point $s$ is known to be $(1-2s)^{n/2}\exp\{s\Delta^2/(1-2s)\}$, we have for all $f,\sigma \in \F_j^k \cap \{\sigma < \sigma_0/2 \}$ such that $d_n(f^t , f) \leq \epsilon_n$  

\begin{multline*}
        P_{f,\sigma}\left[ \frac{\sigma^2}{\sigma_0^2} \sum_{i=1}^n \left\{ \frac{Y_i - f^t_i}{\sigma}\right\}^2 \geq c_2 n  \right]
     \\ \leq \exp\left[ \frac{n}{2} \left\{ -\log(1-2s) + \frac{1}{\sigma^2}\frac{2s}{1-2s}d_n(f,f^t)^2 - 2sc_2 \frac{\sigma_0^2}{\sigma^2}  \right\}  \right] . 
\end{multline*}
For $s$ small enough we have 
\[
\frac{2s}{1-2s} d_n(f,f^t)^2 \leq 4s d_n(f,f^t)^2 \leq 4s \epsilon_n^2 \leq 2sc_2 \sigma_0^2.
\]
Which in turns gives for $c'_2>0$ a fixed constant
\[
\E_{f,\sigma}(1-\Psi_3^t) \leq e^{ -nc_2'} .
\]
Taking $\Psi_3 = \max_t \Psi_3^t$ we get a test such that
\[ 
\E_0^n(\Psi_3) = o(1) ; \sup_{\F_n^j \cap \{ \sigma \leq \sigma_0/2 \} } \E_{f,\sigma} (1-\Psi_3) \leq e^{-Cj^2 n \epsilon_n^2} .
\]
We conclude the proof by taking $\Phi_n= \max\{ \Psi_1, \Psi_2, \Psi_3\}$

\end{proof}

\section{Concentration rate of the posterior distribution for H\"olderian smooth and monotone functions functions}
\label{sec:contraction:rate}

In this section we prove that the posterior concentrate around $f_0, \sigma_0$
at the rate $(n / \log (n))^{- 1 / 4}$ if $f_0 \in \Fm$ and $(n / \log (n))^{-
\alpha / (2 \alpha + 1)}$ if $f_0 \in \mathcal{H}(\alpha, L)$. 

To do so we follow the approach of
{\citet{ghosal:vaart:2007}}. 
Throughout the proof, $C$ will denote a generic constant.

Let $KL(f,g) = \int f \log(f/g)$ be the Kullback-Leibler divergence between the two probability densities $f$ and $g$. We define $V(f,g) = \int (\log(f/g) - KL(f,g))^2 f$. We denote $p_i(\omega,\sigma,k)$ the probability density with respect to the Lebesgue measure of $Y_i = f_{\omega,k} + \epsilon_i$ when $\epsilon_i \sim \mathcal{N}(0,\sigma^2)$ and $p_{i,0}$ the true density of $Y_i$, i.e. when $f = f_0$.
We only consider the case where $f \in \Fm$, a similar proof holds when $f \in \mathcal{H}(\alpha,L)$. We define
\[ B_n (\epsilon) = \left\{ \sum_{i = 1}^n KL\{ p_i(\omega,\sigma,k\}, p_{i,0} ) \leq n \epsilon^2, \sum_{i = 1}^n V\{ p_i(\omega,\sigma,k), p_{i,0} \} \leq n \epsilon^2 \right\} 
\]
Here $p (\omega, \sigma, k)$ and $p_0$ are Gaussian distributions, we can
easily compute

\begin{align*}
KL\{ p_i(\omega,\sigma,k), p_{i,0} \} &= \frac{1}{2} \log\left( \frac{\sigma^2}{\sigma_0^2} \right) - \frac{1}{2}\left( 1 - \frac{\sigma_0^2}{\sigma^2} \right) + \frac{1}{2}
{\left\{f_{\omega,k}(x_i) - f_0(x_i)\right\}^2 \over \sigma^2} \\ 
V\{ p_i(\omega,\sigma,k), p_{i,0} \} &= \frac{1}{2} \left( 1 - \frac{\sigma_0^2}{\sigma^2} \right)^2 + \left[ \frac{\sigma_0^2}{\sigma^2} \left\{ f_{\omega,k} (x_i) - f_0(x_i) \right\}
\right]^2 
\end{align*}

We have $B_n (\epsilon_n) \supset \{d_n^2(f_{\omega, k} , f_0) \leq C
\epsilon_n^2, | \sigma^2 - \sigma_0^2 |^2 \leq C \epsilon_n^2 \}$.\\

For $f_0 \in \F$, denoting $\omega^0_j = n_j^{- 1}  \sum_{x_i \in I_j} f_0
(x_i)$ and $\underline{x_j} = \inf (I_j)$, $\overline{x_j} = \sup (I_j)$ we
have
\[ d_n^2(f_{\omega, k} , f_0 ) = d_n^2( f_0 , f_{\omega^0, k} ) + d_n^2(
   f_{\omega, k} , f_{\omega^0, k}) \]
and

\begin{align*}
d_n^2(f_0 , f_{\omega^0,k})&= \frac{1}{n} \sum_{j=1}^k \sum_{x_i \in I_j} \{f_0(x_i) - f_{\omega^0,k}\}^2 \\
    & \leq \frac{1}{n} \sum_{j=1}^k n_j \{f_0(\underline{x_j}) - f_0(\overline{x_j})\}^2 \\
    & \leq \frac{C}{k} \left[ \sum_{j=1}^{k} \{f_0(\underline{x_j}) -
f_0(\overline{x_j})\} \right]^2 \leq \frac{C ||f_0||_\infty^2}{k} . 
\end{align*} 
Denoting $k_n = C  \lceil ||f_0||^2_\infty \{n / \log (n)\}^{1 / 2} \rceil$ we deduce that $B_n
(\epsilon_n) \supset \{k = k_n, || \omega - \omega^0 ||_{k_n}^2 \leq \epsilon_n^2, |
\sigma^2 - \sigma_0^2 | \leq \epsilon_n^2 \}$ where $|| \cdot ||_{k}$ is the
standard Euclidean norm in $\R^{k}$ i.e. for $a = (a_1,\dots,a_k) \in \R^k$ 
$$
||a||_k^2 = k^{-1} \sum_{i=1}^k a_i^2.
$$
We deduce that for a fixed positive constant $C_0$ that depends on $||f_0||_\infty$ , 
\begin{equation}
  \pi \{B_n (\epsilon_n)\} \gtrsim \left( C \inf_{x \in [0, 1]}
  [g\{f_0 (x)\}] \epsilon_n \right)^{k_n} \pi_{\sigma} (\sigma_0^2)
  \epsilon_n^2 \pi(k = k_n) \geq e^{- C_0 n \epsilon_n^2} \label{eq:piB1} . 
\end{equation}

To end the proof the standard approach of \citet{ghosal:vaart:2007} requires the existence of an exponentially consistent sequence of tests. Their Theorem 4 suited for independent observations relies on the fact that the set $\left\{ d_n(f_{\omega,k} , f_0)^2 + (\sigma - \sigma_0)^2 \geq \epsilon_n^2 \right\} $ can be covered with Hellinger balls. Because of the unknown variance, this cannot be done here, we thus use an alternative approach and to construct tests, and then apply Theorem 3 from \citet{ghosal:vaart:2007}. To prove the existence of tests we apply Lemma \ref{lem:tests}. Note that we have 
\begin{equation}
\F^k_j \subset \left\{ ||\omega - \omega^0||_k \leq Cj\epsilon_n, |\sigma
 - \sigma_0| \leq Cj\epsilon_n \right\},
 \label{eq:Fkj}
\end{equation}
and thus for all $\xi >0$ there exist a $ \xi \epsilon_n$ net of $\F^k_j$ containing less than $(Cj/\xi)^k$ points.

\section{Proof of Lemma 2}

Let $f_0$ either belong to $\F$ or to $\mathcal{H}(\alpha,L)$ and $\epsilon_n$ represent either $\epsilon_n(\F)$ if $f_0 \in \F$ or $\epsilon_n(\alpha)$ if $f_0 \in \mathcal{H}(\alpha,L)$. We denote $A_n = \{(\omega,\sigma,k), d_n(f_{\omega,k},f_0)^2 + |\sigma - \sigma_0|^2 \leq \epsilon_n^2\}$ with $\epsilon_n$ as in Section \ref{sec:contraction:rate}. Thus $\pi(A_n^c | Y_n) = o_{P_0^n}(1)$. We now derive an upper bound for 
$
\pi(\max_j|\omega_j - \omega_j^0| \geq A \xi_n^k |Y_n,A_n).
$ To do so, we look at the following decomposition for all $k_n \in \N$, 
\begin{multline}
  \pi(\max_j|\omega_j - \omega_j^0| \geq A \xi_n^k |Y_n,A_n) \leq  \\ \sum_{k\leq k_n} \pi(k|Y_n,A_n ) \sum_{j=1}^k \int \pi(|\omega_j - \omega_j^0| \geq C \xi_n^k |Y_n,A_n,k,\sigma) d\pi(\sigma|Y_n,A_n,k) + \pi(k>k_n|Y_n).
\label{eq:decompo_2}
\end{multline}
Given Lemma C3 we have, choosing $k_n = C_1 n\epsilon_n^2 $ a constant $C_1$ as in Lemma C3,
$$
\pi(k>k_n|Y_n) = o_{P_0^n}(1)
$$
We now find an upper bound uniformly in $\sigma$ over $A_n$ for $\pi(|\omega_j - \omega_j^0| \geq A \xi_n^k |Y_n,A_n,k,\sigma)$. We first denote $I_l(\omega^0_j,\sigma_0) = \{ l \sigma_0 \xi_n^k \leq | \omega_j - \omega^0_j | \leq (l+1) \sigma_0 \xi_n^k  \} $. We have for $l_0 \leq A$
$$
\Pi(|\omega_j - \omega_j^0| \geq A \xi_n^k |Y_n,A_n,k,\sigma) \leq \sum_{l \geq l_0} \Pi\{I_l(\omega^0_j, \sigma_0) |Y_n,A_n,k,\sigma \}.
$$
We then write
$$
\Pi\{I_l(\omega^0_j, \sigma_0) |Y_n,A_n,k,\sigma\} = { \int_{I_l(\omega^0_j, \sigma_0)} e^{l_n^\sigma(\omega) - l_n^{\sigma_0}(\omega^0) } d\Pi(\omega) \over\int e^{l_n^\sigma(\omega) - l_n^{\sigma_0}(\omega^0)} d\Pi(\omega) },
$$
where $l_n^\sigma(\omega) = -n\log(\sigma^2)/2 - \frac{1}{2} \sum_{i=1}^n \{Y_i - f_{\omega,k}(x_i)\}^2/\sigma^2$. Standard algebra leads to 
$$
l_n^\sigma(\omega) - l_n^{\sigma_0}(\omega^0) = - \frac{1}{2} \sum_{j=1}^k \frac{(\omega_j - \omega_j^0)^2}{\sigma^2} + \sum_{x_i \in I_j}\epsilon_i\frac{\sigma_0}{\sigma^2} (\omega_j - \omega_j^0) + \Delta(\epsilon,\sigma,f_0,k),
$$
where $\Delta(\epsilon,\sigma,f_0,k)$ does not depend on $\omega$ and $\epsilon_i \overset{iid}{\sim} \mathcal{N}(0,1)$ under $p_0^{n}$. We thus deduce 
\begin{multline*}
\Pi\{I_l(\omega^0_j, \sigma_0) |Y_n,A_n,k,\sigma\} =  \\ 
{ \int_{I_l(\omega^0_j, \sigma_0) } \exp\left\{- \frac{1}{2} n_j\frac{(\omega_j- \omega_j^0)^2}{\sigma^2} + \sum_{x_i \in I_j}(\epsilon_i)\frac{\sigma_0}{\sigma^2} (\omega_j - \omega_j^0)\right\} d\Pi(\omega) 
\over
\int \exp\left\{- \frac{1}{2} n_j\frac{(\omega_j- \omega_j^0)^2}{\sigma^2} + \sum_{x_i \in I_j}(\epsilon_i)\frac{\sigma_0}{\sigma^2} (\omega_j - \omega_j^0) \right\} d\Pi(\omega) }  
= \frac{N_{n,j,l}^k(\sigma)}{D_{n,j}^k(\sigma)}
\end{multline*}
We now prove that on a set $\mathcal{E}$ such that $P_0^n(\mathcal{E}) = 1 + o(1)$ we have for $(\epsilon_i) \in \mathcal{E}$, We have an upper bound for $N_{n,j}^k/D_{n,j}^k$ uniformly in $\sigma \in A_n$ for all $k \leq k_n$. 
Let $\mathcal{E} = \left\{ \cap_{k\leq k_n} \cap_{j=1}^k \left\{ \left|\sum_{x_i \in I_j} \epsilon_i \right| \leq c_e \sqrt{n_j\log(n)} \right\} \right\}$ for some constant absolute constant $c_e$ large enough. We compute 
\[
  \mathrm{pr}(\mathcal{E}^c) \leq 2\sum_{k=2}^{k_n} \sum_{j=1}^k \mathrm{pr}\left(\sum_{x_i \in I_j} \epsilon_i > c_e \sqrt{n_j\log(n)}\right) \leq 2\frac{k_n^2}{n^{c_e^2}} = o(1) .
\]

For $(\epsilon_i) \in \mathcal{E}$ and uniformly in $\sigma$ over $A_n$ we compute 
\begin{align*} 
  D_{n,j}^k (\sigma)&= \int\exp \left\{-\frac{n_j}{2\sigma^2}(\omega_j - \omega_j^0)^2 + \frac{\sigma_0}{\sigma^2}(\omega_j - \omega_j^0) \sum_{x_i \in I_j} \epsilon_i  \right\} d \pi(\omega_j) \\ 
  &\geq \int_{|\omega_j - \omega_j^0| \leq \sigma_0 c_e \xi_n^k} \exp \left\{ -n_j(\omega_j - \omega_j^0)^2 - 2 c_e \frac{\sigma_0}{\sigma^2}  n_j |\omega_j - \omega_j^0| \sqrt{\frac{\log(n)}{n_j}} \right\} d \pi(\omega_j) \\ 
  & \geq e^{-3 c_e^2 \sigma_0^2 n_j (\xi_n^k)^2 /(2 \sigma^2) } \Pi(|\omega_j - \omega_j^0| \\ 
 & \geq e^{-3 c_e^2 \sigma_0^2 n_j (\xi_n^k)^2 /(2 \sigma^2) } \sigma_0 c_e \xi_n^k \inf_{x \in [-K,K]}\{g(x)\}  
\end{align*}

Similarly for $(\epsilon_i) \in \mathcal{E}$ and uniformly in $\sigma$ over $A_n$ we have for $l$ large enough 
\begin{align*}
  N_{n,j,l}^k(\sigma) &\leq 
    \int_{I_l(\omega^0_j,\sigma_0)} \exp\left\{- \frac{1}{2} n_j|\omega_j - \omega_j^0| \left( \frac{|\omega_j- \omega_j^0|}{\sigma^2} - \frac{\sigma_0}{\sigma^2} c_e\sqrt{\frac{\log(n)}{n_j}} \right) \right\} d\pi(\omega) \\ 
    &\leq e^{-l^2 \sigma_0^2 n_j (\xi_n^k)^2/(4 \sigma^2) } \Pi\{I_l(\omega_j^0,\sigma_0)\} \\  
     & \leq e^{-l^2 \sigma_0^2 n_j (\xi_n^k)^2/(4 \sigma^2) } \sigma_0 \xi_n^k ||g||_\infty .
\end{align*}

We thus have for $(\epsilon_i)_i \in \mathcal{E}$, $\epsilon >0$ and $l$ large enough, together with condition \textbf{C2} 
\begin{align*} 
\frac{N_{n,j,l}^k(\sigma)}{D_{n,j}^k(\sigma)} &\leq e^{-\frac{1}{2\sigma^2}\sigma_0^2 n_j (\xi_n^k)^2 ( l/2 - 3c_e) } \frac{\Pi\{I_l(\omega_j^0,\sigma_0)\}}{\Pi(|\omega_j - \omega_j^0| \leq \sigma_0 c_e \xi_n^k)  }   \\
&\leq e^{-n_j (\xi_n^k)^2 l^2 \frac{\sigma_0^2}{8\sigma^2} } \frac{||g||_\infty}{c_e \inf_{x \in [-K,K]}\{g(x)\}}, 
\end{align*}
which in turns gives an upper bound for $\Pi(|\omega_j - \omega_j^0| \geq A \xi_n^k |Y_n,A_n,k,\sigma)$ 

$$
\Pi(|\omega_j - \omega_j^0| \geq A \xi_n^k |Y_n,A_n,k,\sigma) \leq \frac{1}{2} e^{- l_0 \frac{\sigma_0^2}{8\sigma^2} n_j (\xi_n^k)^2 } \frac{||g||_\infty}{c_e \inf_{x \in [-K,K]}\{g(x)\}} . 
$$

We thus deduce for $C>0$ an absolute constant and $C'$ depending only on $\pi$,  
$$ \Pi(\max_{1\leq j \leq k} |\omega_j - \omega_j^0| \geq A \xi_n^k |Y_n) \leq C' k_n e^{- l_0 C \log(n) } + o_{P_0^n}(1),
$$
which gives choosing $A$ large enough
$$
P_0^n\left\{ \Pi(\max_{1\leq j \leq k} |\omega_j - \omega_j^0| \geq A \xi_n^k |Y_n) <  {\gamma_1 \over \gamma_0 + \gamma_1} \right\} \to 1. 
$$

\section{Proof of Lemma 3}

Let be either $k_n = n\epsilon_n^2/\log(n)$ if $L(k) = \log(k)$ or $k_n = n\epsilon_n^2$ if $L(k) = 1$. Similarly to before, we have $\pi \left\{ B_n (\epsilon_n) \right\} \geq
e^{- n \epsilon_n^2}$. We define $N_n$ and $D_n$ such that
\[ \pi ( \mathcal{K}_n^c |Y_n) = \frac{\sum_{k \in \mathcal{K}_n^c} \pi (k) 
   \int \frac{p (\omega, \sigma, k)}{p_0} (Y^n) d \Pi (\omega, \sigma)}{\sum_k
   \pi (k)  \int \frac{p (\omega, \sigma, k)}{p_0} (Y^n) d \Pi (\omega,
   \sigma)} = \frac{N_n}{D_n} \]
Given Lemma 10 of {\citet{ghosal:vaart:2007}}, we have
\[ P_0^n  \left( D_n \leq e^{- Cn \epsilon_n^2} \right) = o (1) \]
Note also that
\[ \E_0^n (N_n) = \sum_{k \in \mathcal{K}_n^c} \pi (k)  \int \int_{\R^n}
   \frac{p (\omega, \sigma, k)}{p_0} (Y^n) p_0 (Y_n) d \Pi (\omega, \sigma) d Y^n =
   \pi (k \leq k_n) \leq ce^{- C_u k_n L (k_n)} \]
Thus for $C$ small enough we have

\begin{align*}
\E_0^n \left\{\Pi\left( k \in \mathcal{K}_n^c | Y^n  \right) \right\}&= \E_0^n\left\{
{N_n \over D_n} \I_{D_n > e^{-Cn\epsilon_n^2} }\right\} + o(1) \\ 
& \leq e^{Cn\epsilon_n^2} c e^{-C_u k_n L(k_n)} + o(1) \\ 
&\leq o(1)  
\end{align*}

\end{document}